\documentclass[a4paper,11pt,twoside]{article}
\usepackage{pst-fill,pst-grad}
\usepackage{textcomp}
\usepackage[english]{babel}
\usepackage{graphicx}
\usepackage{amsmath}
\usepackage{float}
\usepackage[matrix,arrow,curve]{xy}
\usepackage{pstricks}
\usepackage{amsmath,amsfonts,verbatim,afterpage,theorem,euscript,mathrsfs,amssymb}
\usepackage{amsfonts}
\usepackage{amssymb}
\usepackage{array}
\usepackage{dsfont}
\usepackage{hyperref}
\def \vu{\vec{u}}
\def \vU{\vec{U}}

\def \vf{\vec{f}}
\def \vg{\vec{g}}
\def \vw{\vec{\omega}}

\def \vn{\vec{\nabla}}
\def \rot{\grad \wedge}
\def \grad{\vn}
\def \R{\mathbb{R}^3}
\newtheorem{Theoreme}{Theorem}

\newtheorem{Proposition}{Proposition}[section]
\newtheorem{Lemme}{Lemma}[section]
\newtheorem{Corollaire}{Corollary}[section]

\numberwithin{equation}{section}
\title{\bf  Some remarks about the stationary Micropolar fluid equations: existence, regularity and uniqueness.} 

\usepackage{authblk}
\author{Diego Chamorro\footnote{\emph{diego.chamorro@univ-evry.fr} (corresponding author)} }
\author{David Llerena\footnote{\emph{david.llerena@univ-evry.fr}} }
\author{Gast\'on Vergara-Hermosilla\footnote{\emph{gaston.vergarahermosilla@univ-evry.fr}} }
\affil{\footnotesize LaMME, Univ. Evry, CNRS, Universit\'e Paris-Saclay, 91025, Evry, France.}
\begin{document}
\maketitle
\begin{scriptsize}
\abstract{We consider here the stationary Micropolar fluid equations which are a particular generalization of the usual Navier-Stokes system where the microrotations of the fluid particles must be taken into account. We thus obtain two coupled equations: one based mainly in the velocity field $\vu$ and the other one based in the microrotation field $\vw$. We will study in this work some problems related to the existence of weak solutions as well as some regularity and uniqueness properties. Our main result establish, under some suitable decay at infinity conditions for the velocity field only, the uniqueness of the trivial solution.}\\[3mm]
\end{scriptsize}
%\tableofcontents
%%%%%%%%%%%%%%%%%%%%%%%%%%%%%%%%%%%%%%%%%%%%%%%%%%%
\section{Introduction and presentation of the results}
In this article we study some problems related to existence, regularity and uniqueness (via Liouville-type theorems) of the 3D stationary Micropolar equations. This system is composed by the usual incompressible stationary 3D Navier-Stokes equations coupled with a stationary equation describing the angular velocity of the rotation of fluid particles.  More precisely, this system of equations reads as follows: 
\begin{equation}\label{SteadyMicropolarEquations}
\begin{cases}
\Delta \vu-(\vu \cdot \vn) \vu-\vn p+\frac{1}{2} \vn \wedge \vw+\vec{f}=0,\qquad div(\vu)=div(\vec{f})=0,\\[5mm]
\Delta \vw+\vn div(\vw)-\kappa\vw-(\vu \cdot \vn) \vw+\frac{1}{2} \vn \wedge \vu
+\vec{g}=0.
\end{cases}
\end{equation}
where $\vu:\R \longrightarrow \R $ denotes the velocity vector field, $p:\R \longrightarrow \mathbb{R} $ denotes the internal pressure of the fluid, $\vw:\R \longrightarrow \R $ is the angular velocity which is usually called the microrotation field and the vector fields $\vf, \vg:\R \longrightarrow \R $ are two given external forces. Here $1\ll \kappa<+\infty$ is a technical parameter. Note that we have the conditions $div(\vu)=div(\vf)=0$ for the first equation above (which is related to the Navier-Stokes system) and this fact implies the following property: since we have $div(\vu)=div(\vn\wedge \vw)=div(\vf)= 0$, formally applying the divergence operator in the first equation we obtain the identity 
\begin{equation}\label{Equation_Pression_intro}
-\Delta p=div\big((\vu \cdot \vn)\vu\big)=div\big(div(\vu\otimes \vu)\big),
\end{equation}
which allows us to recover some information over $p$ from the velocity field $\vu$ independently from $\vw$.\\

Remark also that, as we do not have these conditions for the vector fields $\vw$ and $\vg$, the study of the second equation of (\ref{SteadyMicropolarEquations}) will require a different treatment as we shall see later on. \\

The Micropolar equations were first introduced in the literature in a paper due to Eringen \cite{eringen1966theory}, and nowadays it is applied in studies related to polymers, muddy fluids, nematic liquid crystals, budly liquids, to mention a few. For more applications of this model see the monography \cite{lukaszewicz1999micropolar}. From a more theoretical point of view, the Micropolar system (the evolution problem or the stationary one) was studied in \cite{chamorro2023crypto}, \cite{LUKASZEWICZ200291}, \cite{jens2021regularity}, \cite{Ortega_Rojas99}, \cite{Rojas97}, \cite{RojasBoldrini98} and \cite{wang2020global} (see also the references there in).\\ 

One interesting feature of this model is the fact that it is possible to perform a separated study of the properties of the variables $\vu$ and $\vw$ by imposing different constraints on each one of these variables. In this article we will follow this direction to study the uniqueness of weak solutions of the system (\ref{SteadyMicropolarEquations}). To do so, we will first establish an existence result and we will study some regularity properties. Indeed, assuming some mild assumptions on the external forces $\vf, \ \vg$, we have the following theorem:
%%%%%%%%%%%%%%%%%%%%%%%%%%%%%%%%%%%%%%%%%%%%%%%%%%%
\begin{Theoreme}[Existence]\label{Theorem_Existence}
Let  $\vec{f},\vec{g}\in \dot H^{-1}(\mathbb{R}^3)\cap \dot H^{-2}(\mathbb{R}^3)$ be two exterior forces (with $div(\vf)=0$ but $div(\vg)\neq 0$). Then, there exist at least one weak solution $(\vu,p,\vw)$ of the stationary Micropolar equations \eqref{SteadyMicropolarEquations} such that we have $\vu\in \dot{H}^1(\mathbb{R}^3)$, $\vw\in H^1(\R)$ and  $p\in  \dot H^{\frac{1}{2}}(\mathbb{R}^3)$.  
\end{Theoreme}
%%%%%%%%%%%%%%%%%%%%%%%%%%%%%%%%%%%%%%%%%%%%%%%%%%%
Existence of certain weak solutions to the stationary micropolar equations \eqref{SteadyMicropolarEquations} has been studied in the literature. For instance,  in \cite{LUKASZEWICZ200291} the authors deal the case of 3D bounded domains with boundary data in $L^2(\R)$ and for exterior domains we refer to \cite{DuOrRo03}. Our approach here is quite different as we study the existence of solutions in the whole space $\R$: we thus start by mollifying the system \eqref{SteadyMicropolarEquations} and we study the existence of solutions using the well known Scheafer's fixed point theorem, then by a suitable limit process we can get rid of the mollification and we obtain a weak solution of the initial problem. Let us stress that we do not claim any optimality on the functional spaces considered here as other more general frameworks can be surely used.\\ 

As we mentioned above, the setting given in Theorem \ref{Theorem_Existence} is quite interesting for our purposes as we can deduce the following regularity result:
%%%%%%%%%%%%%%%%%%%%%%%%%%%%%%%%%%%%%%%%%%%%%%%%%%%
\begin{Theoreme}[Regularity]\label{Theorem_Regularity}
Assume that $\vf=\vg=0$ or that $\vf$ and $\vg$ are regular vector fields. Let $(\vu,p,\vw)$ be a solution of the stationary Micropolar equations \eqref{SteadyMicropolarEquations} obtained in the Theorem  \ref{Theorem_Existence} above. Then the functions $\vu$,  $p$ and $\vw$ are regular. 
\end{Theoreme}
%%%%%%%%%%%%%%%%%%%%%%%%%%%%%%%%%%%%%%%%%%%%%%%%%%%
The regularity of the system will be obtained by an iterative process: we will start by an independent study of each variable, but as it is a coupled system, at some point the regularity of each variable will depend on the regularity of the other variables. Thus, by a suitable bootstrap argument we will obtain the wished gain of regularity.\\

Let us now consider the problem of uniqueness of weak solutions of the Micropolar equations. In the case when $\vf=\vg=0$, it is easy to see that $\vu=0$, $p=0$ and $\vw=0$ is a solution of the system \eqref{SteadyMicropolarEquations}. But this trivial solution is not unique: indeed, if we define the function $\psi:\R \longrightarrow \mathbb{R}$ by $\psi(x_1,x_2,x_3)=\frac{x^{2}_{1}}{2}+\frac{x^{2}_{2}}{2}-x^{2}_{3}$ and if we set the functions $\vu$, $p$ and $\vw$ by the identities
$$\vu(x_1,x_2,x_3)=\vn\psi(x_1,x_2,x_3)=(x_1, x_2, -2x_3),\quad p(x_1,x_2,x_3)=-\frac{1}{2}\| \vu(x_1,x_2,x_3)\|^2 \quad \mbox{and}\quad \vw=0,$$
then, using basic rules of vector calculus we have that $(\vu, p, \vw)$ given by the expressions above satisfies in the weak sense the system (\ref{SteadyMicropolarEquations}).\\

Note furthermore that even in the setting given by the Theorems \ref{Theorem_Existence} and \ref{Theorem_Regularity} (where we have some decay at infinity for the vector fields $\vu$ and $\vw$ and regularity properties for all the variables), if we set $\vf=\vg=0$, we can not ensure that the trivial solution is unique. This fact is of course related to the uniqueness of the stationary Navier-Stokes equation (and related systems) which has been studied intensively in \cite{chae2014liouville}, \cite{chae2013liouville}, \cite{chamorro2021some}, \cite{jarrin2020remark}, \cite{koch2009liouville}, \cite{seregin2016liouville}and \cite{yuan2020liouville} and it is a very challenging open problem.\\

To overcome this issue, and following the pionering work of Galdi \cite{galdi2011introduction} for the stationary Navier-Stokes equations, we will consider an additional condition in order to obtain that the trivial solution is unique when $\vf=\vg=0$. Indeed, if we assume some extra information (stated in terms of a suitable Banach space $E$ that will provide a nice decay at infinity) we shall prove a uniqueness result.  These kind of results are known in the literature as \emph{Liouville-type} theorems.\\ 

We can of course assume some information on the variables $\vu$ and $\vw$, but the main novelty of this article is to assume an extra hypothesis for \emph{only one variable} (namely the velocity field $\vu$) and with this single condition we will see how to deduce that the trivial solution is unique: 
%%%%%%%%%%%%%%%%%%%%%%%%%%%%%%%%%%%%%%%%%%%%%%%%%%%
\begin{Theoreme}[Uniqueness]\label{Theorem_Liouville_type} Set $\vf=\vg=0$ and consider $(\vu,p,\vw)$ a solution of the stationary Micropolar fluids equations \eqref{SteadyMicropolarEquations}
obtained via the Theorem \ref{Theorem_Existence}. Assume in addition that we have the condition $\vu\in L^q(\R)$ with $3\le q\le \frac{9}{2}$. Then we have $\vu=\vw=0$.
\end{Theoreme}
%%%%%%%%%%%%%%%%%%%%%%%%%%%%%%%%%%%%%%%%%%%%%%%%%%%
Let us note that, as we are in the framework of the Theorem \ref{Theorem_Existence}, we have $\vu\in \dot{H}^1(\R)$ and by the Sobolev embeddings we also have that $\vu\in L^6(\R)$. However, this information is not enough to obtain the wished result and we need to impose a stronger decay at infinity with the condition $\vu\in L^q(\R)$ with $3\le q\le \frac{9}{2}$. Following what is known in the Navier-Stokes equations, the upper bound $q\le \frac{9}{2}$ seems to be the best available up to now and the case $\frac{9}{2}<q<6$ is a completely open problem for both systems. Remark also that no particular condition is asked for the variable $\vw$ and this separated study of the variables is, to the best of our knowledge, new for this type of problem for equations (\ref{SteadyMicropolarEquations}). Finally, let us stress that we do not claim any optimality of our result: the additional condition stated here in terms of Lebesgue spaces can probably be replaced by other functional spaces such as Lorentz, Besov or Morrey spaces, etc.\\

The plan of the article is simple: in Section \ref{Secc_Theo1} we prove Theorem \ref{Theorem_Existence}, Section \ref{Secc_Theo2} is devoted to the proof of the Theorem \ref{Theorem_Regularity} and, finally, in Section \ref{Secc_Theo3} we study the Liouville-type result stated in the Theorem \ref{Theorem_Liouville_type}.
%%%%%%%%%%%%%%%%%%%%%%%%%%%%%%%%%%%%%%%%%%%%%%%%%%%
\section{Proof of Theorem \ref{Theorem_Existence}}\label{Secc_Theo1}
Recall that for a regular enough vector field $\vec{\phi}$ we may define the Leray projector by the formula $\mathbb{P}(\vec{\phi})=\vec{\phi}+\vn \frac{1}{(-\Delta)}(\vn\cdot \vec{\phi})$. Among its properties we have the identity $\mathbb{P}(\vn\wedge \vec{\phi})=0$ and if $\vec{\phi}$ is a divergence free vector field ($div(\vec{\phi})=0$) we also have $\mathbb{P}(\vec{\phi})=\vec{\phi}$. Finally if $\varphi:\R\longrightarrow \mathbb{R}$ is a real function, we have 
$\mathbb{P}(\vn \varphi)=0$. With these properties in mind, and since $div(\vu)=div(\vf)=0$, we can apply the Leray projector to the first equation of \eqref{SteadyMicropolarEquations} to obtain
\begin{equation}\label{equationwithoutpressure}
\begin{cases}
0=  \Delta \vu-\mathbb{P}[(\vu \cdot \vn) \vu]+\frac{1}{2}\rot\vw+\vf, 
\\[5mm]
0= \Delta \vw+\vn div(\vw)-(\vu \cdot \grad) \vw-\kappa\vw+\frac{1}{2} \rot \vu+g.
\end{cases}
\end{equation}
Note that the pressure $p$ is absent in the previous system but using the equation (\ref{Equation_Pression_intro}) we obtain the formula
\begin{equation*}
p=\frac{1}{(-\Delta)}\left(div(div (\vu\otimes\vu))\right),
\end{equation*}
which will be helpful to recover the pressure from the velocity $\vu$.\\ 

Now, for $R>1$ we define the auxiliary function $\theta_R(x)=\theta (\frac{x}{R})$ where $\theta \in \mathcal{C}_0^{\infty} (\R)$ is a regular function such that $0\leq \theta (x) \leq 1$, $\theta (x)=1$ for $|x|\leq 1$ and $\theta (x)=0$ for $|x|>2$. Note in particular that we have $\|\theta_R\|_{L^\infty}\leq 1$.\\

For a fixed $0<\epsilon\leq 1$, we then consider the system
\begin{equation}\label{Probleme_Approx}
\begin{cases}
0=-\epsilon \Delta^2 \vu + \Delta \vu-\mathbb{P}\big[\big((\theta_R\vu) \cdot \vn\big) (\theta_R\vu)\big]+\frac{1}{2} \theta^2_R(\vn \wedge \vw)+\vf,
\\[5mm]
0=-\epsilon \Delta^2 \vw + \Delta \vw+\vn (\theta_Rdiv(\vw))-\big((\theta_R\vu) \cdot \vn\big) (\theta_R\vw)-\kappa(\theta^2_R\vw)+\frac{1}{2}\theta^2_R (\vn \wedge \vu)+\vg,
\end{cases}
\end{equation}
of course if we let $\epsilon\to 0$ and $R\to+\infty$ we recover (at least formally) the initial problem (\ref{equationwithoutpressure}). It is easy to see that, if we define $\vU= \begin{pmatrix} \vu \\ \vw \end{pmatrix}$, then the system above can be rewritten in the following form
\begin{equation}\label{FixedPoint}
\vU = T_{R,\epsilon} (\vU),
\end{equation}
where 
\begin{equation}\label{eq_epsilon}
T_{R,\epsilon} (\vU) = 
\dfrac{1}{[\epsilon \Delta^2  + (-\Delta) ]} 
\left(\begin{array}{c}  
-\mathbb{P}\left[\left(\left(\theta_R \vu\right) \cdot \vn\right)\left(\theta_R \vu\right)\right]
+ \frac{1}{2}\theta^2_R (\vn \wedge \vw)+\vf\\ [5mm]
\vn (\theta_R  div(\vw))-\left(\left(\theta_R \vu\right) \cdot \vn\right)\left(\theta_R \vw\right)-\kappa(\theta^2_R \vw)+\frac{1}{2}\theta^2_R (\vn \wedge \vu)+  \vec{g}
\end{array}\right).
\end{equation}
Hence, our strategy consists in  applying  a fixed point theorem to the operator $T_{\epsilon,R}$ in order to obtain a solution of the mollified problem (\ref{FixedPoint}). For this, we will use the following result, which is a variant of the Leray-Schauder theorem (see Theorem 16.1 in \cite{lemarie2018navier}):
%%%%%%%%%%%%%%%%%%%%%%%%%%%%%%%%%%%%%%%%%%%%%%%%%%%
\begin{Theoreme}[Schaefer]\label{Teo_SchauderFixPoint}
Let $(E,\|\cdot\|_E)$ be a Banach space and $T:E\longrightarrow E$ an application which is continuous and compact. If there exists a finite $M$ such that, for every $\lambda \in[0,1]$, $e = \lambda T(e)$ implies $\|e\|_E\leq M$, then there exists at least one $e \in E$ such that $T (e) = e$.
\end{Theoreme}
%%%%%%%%%%%%%%%%%%%%%%%%%%%%%%%%%%%%%%%%%%%%%%%%%%%
Therefore, in order to obtain a solution for the problem (\ref{FixedPoint}), we only need to prove that the operator $T_{R,\epsilon}$ given in \eqref{eq_epsilon} verifies the hypothesis of the Theoreme \ref{Teo_SchauderFixPoint} above. For this purpose, we will divide the study in three propositions: first we study the continuity of the application $T_{R,\epsilon}$, next we investigate its compactness and, finally, we shall prove the a priori estimates that will allow us to obtain the wished result. 
%%%%%%%%%%%%%%%%%%%%%%%%%%%%%%%%%%%%%%%%%%%%%%%%%%%
\begin{Proposition}[Continuity]\label{Propo_Continuity}
Under the general setting of the Theorem \ref{Theorem_Existence}, the application $T_{R,\epsilon}$ defined in \eqref{eq_epsilon}  is continuous in the space $\left(\dot H^1(\R)\cap \dot H^2(\R) \right)^2$.
\end{Proposition}
%%%%%%%%%%%%%%%%%%%%%%%%%%%%%%%%%%%%%%%%%%%%%%%%%%%
{\bf Proof.} Let us consider the space $E=\dot{H}^1(\mathbb{R}^3)\cap \dot{H}^2(\mathbb{R}^3)$, endowed with the norm 
$$\|\cdot\|_E= \|\cdot\|_{\dot{H}^1} + \sqrt{\epsilon}\|\cdot\|_{\dot{H}^2},$$
 and we define the subspace $E_\sigma$ of $E$ by 
$$E_\sigma=\big\{\vec{\phi}:\R\longrightarrow\R: \vec{\phi}\in E, \; div(\vec{\phi})=0\big\},$$
this space will be  endowed with the same norm $\|\cdot\|_{E}$, but for the sake of clarity will be denoted it by $\|\cdot\|_{E_\sigma}$. We will now study the continuity of the application $T_{R,\epsilon} $ in the functional space $E_\sigma\times E$, normed by $\|(\vu,\vw) \|_{E_\sigma\times E} = \|\vu\|_{E_\sigma} + \|\vw\|_{E}$. We write then
\begin{eqnarray*}
&&\|T_{R,\epsilon}(\vU)\|_{E_\sigma\times E}=\\
&&\left\|\dfrac{1}{[\epsilon \Delta^2  + (-\Delta) ]} 
\left(\begin{array}{c}  
-\mathbb{P}\left[\left(\left(\theta_R \vu\right) \cdot \vn\right)\left(\theta_R \vu\right)\right]
+ \frac{1}{2}\theta^2_R (\vn \wedge \vw)+\vf\\ [5mm]
\vn (\theta_R  div(\vw))-\left(\left(\theta_R \vu\right) \cdot \vn\right)\left(\theta_R \vw\right)-\kappa(\theta^2_R \vw)+\frac{1}{2}\theta^2_R (\vn \wedge \vu)+  \vec{g}
\end{array}\right)\right\|_{E_\sigma\times E},
\end{eqnarray*}
and we obtain
\begin{eqnarray}
&\leq &\underbrace{\left\|\dfrac{1}{[\epsilon \Delta^2+(-\Delta) ]} \mathbb{P}\left[\left(\left(\theta_R \vu\right) \cdot \vn\right)\left(\theta_R \vu\right)\right]\right\|_{E_\sigma}}_{(1)}+\frac12\underbrace{\left\|\dfrac{1}{[\epsilon \Delta^2+(-\Delta) ]} \theta^2_R (\vn \wedge \vw)\right\|_{E_\sigma}}_{(2)}+\underbrace{\left\|\dfrac{1}{[\epsilon \Delta^2+(-\Delta) ]} \vf\right\|_{E_\sigma}}_{(3)}\notag\\
&+&\underbrace{\left\|\dfrac{1}{[\epsilon \Delta^2+(-\Delta) ]}\vn (\theta_R  div(\vw))\right\|_{E}}_{(4)}+\underbrace{\left\|\dfrac{1}{[\epsilon \Delta^2+(-\Delta) ]}\left(\left(\theta_R \vu\right)\cdot \vn\right)\left(\theta_R \vw\right)\right\|_{E}}_{(5)}+\underbrace{\left\|\dfrac{1}{[\epsilon\Delta^2+(-\Delta) ]}\kappa(\theta^2_R \vw)\right\|_{E}}_{(6)}\notag\\
&+&\frac{1}{2}\underbrace{\left\|\dfrac{1}{[\epsilon \Delta^2+(-\Delta) ]} \theta^2_R (\vn \wedge \vu)\right\|_{E}}_{(7)}+\underbrace{\left\|\dfrac{1}{[\epsilon \Delta^2+(-\Delta) ]} \vg\right\|_{E}}_{(8)}.\label{Estimation_Continuite01}
\end{eqnarray}
Before moving on to the study of each of these terms we will state a lemma that will be useful in the sequel:
%%%%%%%%%%%%%%%%%%%%%%%%%%%%%%%%%%%%%%%%%%%%%%%%%%%
\begin{Lemme}\label{Lemme_EstimationFourier}
Let $0<\epsilon<1$ and consider a parameter $1\leq \sigma \leq 2$. Then, for all function $\vec{\phi}\in L^2(\R)$ we have the estimates 
$$ \left\|\dfrac{(-\Delta)^\sigma}{[\epsilon \Delta^2+(-\Delta)]}\vec{\phi}\right\|_{L^2}\leq \frac{C}{\epsilon}\|\vec{\phi}\|_{L^2},$$
for some constant $C>1$, in particular if $\sigma=1$ the constant does not depend on $\epsilon$.
\end{Lemme}
%%%%%%%%%%%%%%%%%%%%%%%%%%%%%%%%%%%%%%%%%%%%%%%%%%%
Indeed, since the Fourier symbol associated to the operator $\dfrac{(-\Delta)^\sigma}{[\epsilon \Delta^2+(-\Delta)]}$ is $\frac{|\xi|^{2\sigma}}{\epsilon|\xi|^4+|\xi|^{2}}$, which is controlled uniformly (since $1\leq \sigma\leq 2$) by $\frac{C}{\epsilon}$ we obtain the announced estimate in the space $L^2(\R)$.\\

With this lemma at hand, we can study now all the terms in (\ref{Estimation_Continuite01}):
\begin{itemize}
%%%%%%%%%%%%%%%%%%%%%%%%%%%%%%%%%%%%%%%%%%%%%%%%%%%
\item For the first term of the expression (\ref{Estimation_Continuite01}),  we have:
\begin{eqnarray*}
\left\|\dfrac{1}{[\epsilon \Delta^2+(-\Delta) ]} \mathbb{P}\left[\left(\left(\theta_R \vu\right) \cdot \vn\right)\left(\theta_R \vu\right)\right]\right\|_{E_\sigma}&=&\left\|\dfrac{1}{[\epsilon \Delta^2+(-\Delta) ]} \mathbb{P}\left[\left(\left(\theta_R \vu\right) \cdot \vn\right)\left(\theta_R \vu\right)\right]\right\|_{\dot{H}^1}\\
&&+\sqrt{\epsilon}\left\|\dfrac{1}{[\epsilon \Delta^2+(-\Delta) ]} \mathbb{P}\left[\left(\left(\theta_R \vu\right) \cdot \vn\right)\left(\theta_R \vu\right)\right]\right\|_{\dot{H}^2}.
\end{eqnarray*}
Now, by the boundedness properties of the Leray projector and introducing the operators $(-\Delta)$ and $(-\Delta)^{\frac{3}{2}}$, we may write
$$\leq\left\|\dfrac{(-\Delta) }{[\epsilon \Delta^2+(-\Delta) ]} \frac{1}{(-\Delta) }\left(\left(\theta_R \vu\right) \cdot \vn\right)\left(\theta_R \vu\right)\right\|_{\dot{H}^1}+\sqrt{\epsilon}\left\|\dfrac{(-\Delta)^{\frac{3}{2}}}{[\epsilon \Delta^2+(-\Delta) ]} \frac{1}{(-\Delta)^{\frac{3}{2}}}\left(\left(\theta_R \vu\right) \cdot \vn\right)\left(\theta_R \vu\right)\right\|_{\dot{H}^2},$$
and applying the Lemma \ref{Lemme_EstimationFourier} with $\sigma$ equal to $1$ and $\frac{3}{2} $ respectively,  we have 
\begin{eqnarray}
&\leq &C\left\|\frac{1}{(-\Delta) }\left(\left(\theta_R \vu\right) \cdot \vn\right)\left(\theta_R \vu\right)\right\|_{\dot{H}^1}+\left(\frac{C}{\epsilon}\right)\sqrt{\epsilon}\left\| \frac{1}{(-\Delta)^{\frac{3}{2}}}\left(\left(\theta_R \vu\right) \cdot \vn\right)\left(\theta_R \vu\right)\right\|_{\dot{H}^2}\notag\\
&\leq & C\left\|\left(\left(\theta_R \vu\right) \cdot \vn\right)\left(\theta_R \vu\right)\right\|_{\dot{H}^{-1}}+\frac{C}{\sqrt{\epsilon}}\left\|\left(\left(\theta_R \vu\right) \cdot \vn\right)\left(\theta_R \vu\right)\right\|_{\dot{H}^{-1}}\notag\\
&\leq & \left(C+\frac{C}{\sqrt{\epsilon}}\right)\left\|\left(\left(\theta_R \vu\right) \cdot \vn\right)\left(\theta_R \vu\right)\right\|_{\dot{H}^{-1}}.\label{Estimate_Term1}
\end{eqnarray}
By the inclusion $L^{\frac{6}{5}} (\R)\subset \dot{H}^{-1}(\R)$ and the H\"older inequality with $\frac{5}{6} = \frac{1}{3}+\frac{1}{2}$, we obtain
\begin{eqnarray*}
\left\|\left[\left(\theta_R \vu\right) \cdot \vn\right]\left(\theta_R \vu\right)\right\|_{\dot{H}^{-1}}
&\leq &\sum_{i=1}^3\left\|\left(\theta_R u_i\right) \partial_{x_i}\left(\theta_R \vu\right)\right\|_{L^{\frac{6}{5}}}\leq  \sum_{i=1}^3\left\|\theta_R u_i\right\|_{L^3}\left\|\partial_{x_i}\left(\theta_R \vu\right)\right\|_{L^2}\\
&\leq & \sum_{i=1}^3\|\theta_R\|_{L^6} \|u_i\|_{L^6}\big(\|\partial_{x_i}\theta_R\|_{L^3} \|\vu\|_{L^6}+\|\theta_R\|_{L^\infty}\|\vu\|_{\dot{H}^1}\big),
\end{eqnarray*}
where in the last line we used the H\"older inequalities with $\frac{1}{3}=\frac{1}{6}+\frac{1}{6}$ and $\frac{1}{2}=\frac{1}{3}+\frac{1}{6}$. Now, by the Sobolev embedding $\dot H^1 (\R)\subset L^6(\R)$, we easily obtain the estimate
$$\left\|\left[\left(\theta_R \vu\right) \cdot \vn\right]\left(\theta_R \vu\right)\right\|_{\dot{H}^{-1}}\leq C_R\|\vu\|_{\dot{H}^1}\|\vu\|_{\dot{H}^1},$$
from which we deduce 
\begin{equation}\label{Estimation_Continuite1}
\left\|\dfrac{1}{[\epsilon \Delta^2+(-\Delta) ]} \mathbb{P}\left[\left(\left(\theta_R \vu\right) \cdot \vn\right)\left(\theta_R \vu\right)\right]\right\|_{E_\sigma}\leq C_{R, \epsilon}\|\vu\|_{\dot{H}^1}\|\vu\|_{\dot{H}^1}\leq  C_{R, \epsilon}\|\vu\|_{E_\sigma}\|\vu\|_{E_\sigma}.
\end{equation}
%%%%%%%%%%%%%%%%%%%%%%%%%%%%%%%%%%%%%%%%%%%%%%%%%%%
\item For the term (2) in (\ref{Estimation_Continuite01}) we have, using the identity 
$\theta^2_R ( \vn \wedge \vw) = \vn \wedge (\theta^2_R \vw)-  \vw \wedge \vn (\theta^2_R)$:
\begin{eqnarray*}
\left\|\dfrac{1}{[\epsilon \Delta^2+(-\Delta) ]} \theta^2_R (\vn \wedge \vw)\right\|_{E_\sigma}\leq \left\|\dfrac{1}{[\epsilon \Delta^2+(-\Delta) ]} \vn \wedge (\theta^2_R \vw)\right\|_{\dot{H}^1}+\left\|\dfrac{1}{[\epsilon \Delta^2+(-\Delta) ]} \vw \wedge \vn (\theta^2_R)\right\|_{\dot{H}^1}\\
+\sqrt{\epsilon}\left\|\dfrac{1}{[\epsilon \Delta^2+(-\Delta) ]} \vn \wedge (\theta^2_R \vw)\right\|_{\dot{H}^2}+\sqrt{\epsilon}\left\|\dfrac{1}{[\epsilon \Delta^2+(-\Delta) ]} \vw \wedge \vn (\theta^2_R)\right\|_{\dot{H}^2},
\end{eqnarray*}
and introducing the operators $(-\Delta)$ and $(-\Delta)^{\frac{3}{2}}$ in the terms above we can write
\begin{eqnarray*}
&\leq &\left\|\dfrac{(-\Delta)}{[\epsilon \Delta^2+(-\Delta) ]}\frac{1}{(-\Delta)} \vn \wedge (\theta^2_R \vw)\right\|_{\dot{H}^1}+\left\|\dfrac{(-\Delta)}{[\epsilon \Delta^2+(-\Delta) ]}\frac{1}{(-\Delta)} \vw \wedge \vn (\theta^2_R)\right\|_{\dot{H}^1}\\
&+&\sqrt{\epsilon}\left\|\dfrac{(-\Delta)^{\frac{3}{2}}}{[\epsilon \Delta^2+(-\Delta) ]}\frac{1}{(-\Delta)^{\frac{3}{2}}} \vn \wedge (\theta^2_R \vw)\right\|_{\dot{H}^2}+\sqrt{\epsilon}\left\|\dfrac{(-\Delta)^{\frac{3}{2}}}{[\epsilon \Delta^2+(-\Delta) ]}\frac{1}{(-\Delta)^{\frac{3}{2}}} \vw \wedge \vn (\theta^2_R)\right\|_{\dot{H}^2},
\end{eqnarray*}
thus, the Lemma \ref{Lemme_EstimationFourier} we obtain
\begin{eqnarray*}
&\leq &C\left\|\frac{1}{(-\Delta)} \vn \wedge (\theta^2_R \vw)\right\|_{\dot{H}^1}+C\left\|\frac{1}{(-\Delta)} \vw \wedge \vn (\theta^2_R)\right\|_{\dot{H}^1}\\
&&+\frac{C}{\epsilon}\sqrt{\epsilon}\left\|\frac{1}{(-\Delta)^{\frac{3}{2}}} \vn \wedge (\theta^2_R \vw)\right\|_{\dot{H}^2}+\frac{C}{\epsilon}\sqrt{\epsilon}\left\|\frac{1}{(-\Delta)^{\frac{3}{2}}} \vw \wedge \vn (\theta^2_R)\right\|_{\dot{H}^2}.
\end{eqnarray*}
By the properties of the homogeneous Sobolev spaces we have the estimate
\begin{eqnarray}
&\leq &C\|\theta^2_R \vw\|_{L^2}+C\| \vw \wedge \vn (\theta^2_R)\|_{\dot{H}^{-1}}+\frac{C}{\sqrt{\epsilon}}\| \theta^2_R \vw\|_{L^2}+\frac{C}{\sqrt{\epsilon}}\| \vw \wedge \vn (\theta^2_R)\|_{\dot{H}^{-1}}\notag\\
&\leq &C\|\theta^2_R \vw\|_{L^2}+C\| \vw \wedge \vn (\theta^2_R)\|_{L^{\frac{6}{5}}}+\frac{C}{\sqrt{\epsilon}}\|  \theta^2_R \vw\|_{L^2}+\frac{C}{\sqrt{\epsilon}}\| \vw \wedge \vn (\theta^2_R)\|_{L^{\frac{6}{5}}},\label{Estimate_Term2}
\end{eqnarray}
where in the last inequality we used the embedding $L^{\frac{6}{5}}(\R)\subset \dot{H}^{-1}(\R)$. Now, by the H\"older inequalities with $\frac{1}{2}=\frac{1}{3}+\frac{1}{6}$ and $\frac{5}{6}=\frac{1}{6}+\frac{2}{3}$, we can write
$$\leq \left(C+\frac{C}{\sqrt{\epsilon}}\right)\|\theta^2_R\|_{L^3}\|\vw\|_{L^6}+ \left(C+\frac{C}{\sqrt{\epsilon}}\right)\| \vw\|_{L^{6}}\|\vn (\theta^2_R)\|_{L^{\frac{3}{2}}},$$
and using the embedding $\dot{H}^1(\R)\subset L^{6}(\R)$ we have
\begin{equation}\label{Estimation_Continuite2}
\left\|\dfrac{1}{[\epsilon \Delta^2+(-\Delta) ]} \theta^2_R (\vn \wedge \vw)\right\|_{E_\sigma}\leq C_{R,\epsilon}\|\vw\|_{\dot{H}^1}\leq C_{R,\epsilon}\|\vw\|_{E}.
\end{equation}
%%%%%%%%%%%%%%%%%%%%%%%%%%%%%%%%%%%%%%%%%%%%%%%%%%%
\item The quantity (3) in (\ref{Estimation_Continuite01}) is treated as follows, we write:
$$\left\|\dfrac{1}{[\epsilon \Delta^2+(-\Delta)]}\vf\right\|_{E_\sigma}=\left\|\dfrac{(-\Delta)}{[\epsilon \Delta^2+(-\Delta) ]}\frac{1}{(-\Delta)}\vf\right\|_{\dot{H}^1}+\sqrt{\epsilon}\left\|\dfrac{\Delta^2}{[\epsilon \Delta^2+(-\Delta) ]} \frac{1}{\Delta^2}\vf\right\|_{\dot{H}^2},$$
and by the Lemma \ref{Lemme_EstimationFourier} we obtain
\begin{equation*}
\leq C\left\|\frac{1}{(-\Delta)}\vf\right\|_{\dot{H}^1}+\frac{C}{\epsilon}\sqrt{\epsilon}\left\|\frac{1}{\Delta^2}\vf\right\|_{\dot{H}^2}\leq C\|\vf\|_{\dot{H}^{-1}}+\frac{C}{\sqrt{\epsilon}}\|\vf\|_{\dot{H}^{-2}}<+\infty,
\end{equation*}
which is bounded by the hypothesis over the external force $\vf$.
%%%%%%%%%%%%%%%%%%%%%%%%%%%%%%%%%%%%%%%%%%%%%%%%%%%
\item For the term (4) of (\ref{Estimation_Continuite01}) we proceed as follows:
\begin{eqnarray*}
\left\|\dfrac{1}{[\epsilon \Delta^2+(-\Delta) ]}\vn (\theta_R  div(\vw))\right\|_{E}&=&\left\|\dfrac{(-\Delta)}{[\epsilon \Delta^2+(-\Delta) ]}\frac{1}{(-\Delta)}\vn (\theta_R  div(\vw))\right\|_{\dot{H}^1}\\
&&+\sqrt{\epsilon}\left\|\dfrac{(-\Delta)^{\frac{3}{2}}}{[\epsilon \Delta^2+(-\Delta) ]}\frac{1}{(-\Delta)^{\frac{3}{2}}}\vn (\theta_R  div(\vw))\right\|_{\dot{H}^2}.
\end{eqnarray*}
Thus, by the Lemma \ref{Lemme_EstimationFourier}, we can write
$$\leq C\left\|\frac{1}{(-\Delta)}\vn (\theta_R  div(\vw))\right\|_{\dot{H}^1}+\frac{C}{\epsilon}\sqrt{\epsilon}\left\|\frac{1}{(-\Delta)^{\frac{3}{2}}}\vn (\theta_R  div(\vw))\right\|_{\dot{H}^2},$$
and using the properties of the Sobolev space we obtain
$$\leq C\|\theta_R  div(\vw)\|_{L^2}+\frac{C}{\sqrt{\epsilon}}\|\theta_R  div(\vw)\|_{L^2}=\left(C+\frac{C}{\sqrt{\epsilon}}\right)\|\theta_R  div(\vw)\|_{L^2}\leq \left(C+\frac{C}{\sqrt{\epsilon}}\right)\|\theta_R\|_{L^\infty}\|\vw\|_{\dot{H}^1},$$
from which we easily deduce that 
$$\left\|\dfrac{1}{[\epsilon \Delta^2+(-\Delta) ]}\vn (\theta_R  div(\vw))\right\|_{E}\leq C_{R,\epsilon}\|\vw\|_E.$$
%%%%%%%%%%%%%%%%%%%%%%%%%%%%%%%%%%%%%%%%%%%%%%%%%%%
\item The term (5) of (\ref{Estimation_Continuite01}) is given by $\left\|\dfrac{1}{[\epsilon \Delta^2+(-\Delta) ]}\left(\left(\theta_R \vu\right)\cdot \vn\right)\left(\theta_R \vw\right)\right\|_{E}$ and it can be treated in the same fashion as the first term of (\ref{Estimation_Continuite01}): indeed, with the same arguments we obtain (see the estimate (\ref{Estimation_Continuite1})):
$$\left\|\dfrac{1}{[\epsilon \Delta^2+(-\Delta) ]}\left(\left(\theta_R \vu\right)\cdot \vn\right)\left(\theta_R \vw\right)\right\|_{E}\leq  C_{R, \epsilon}\|\vu\|_{E_\sigma}\|\vw\|_{E}.$$
%%%%%%%%%%%%%%%%%%%%%%%%%%%%%%%%%%%%%%%%%%%%%%%%%%%
\item For the term (6) of (\ref{Estimation_Continuite01}) we write
$$\left\|\dfrac{1}{[\epsilon\Delta^2+(-\Delta) ]}\kappa(\theta^2_R \vw)\right\|_{E}=\left\|\dfrac{(-\Delta)}{[\epsilon\Delta^2+(-\Delta)]}\frac{1}{(-\Delta)}\kappa(\theta^2_R \vw)\right\|_{\dot{H}^1}+\sqrt{\epsilon}\left\|\dfrac{(-\Delta)^{\frac{3}{2}}}{[\epsilon\Delta^2+(-\Delta)]}\frac{1}{(-\Delta)^{\frac{3}{2}}}\kappa(\theta^2_R \vw)\right\|_{\dot{H}^2},$$
and applying once more the Lemma \ref{Lemme_EstimationFourier} we obtain
\begin{equation}\label{Estimation_Continuite_Term6}
\leq C\left\|\frac{1}{(-\Delta)}\kappa(\theta^2_R \vw)\right\|_{\dot{H}^1}+\frac{C}{\epsilon}\sqrt{\epsilon}\left\|\frac{1}{(-\Delta)^{\frac{3}{2}}}\kappa(\theta^2_R \vw)\right\|_{\dot{H}^2}\leq \kappa\left(C+\frac{C}{\sqrt{\epsilon}}\right)\left\|\theta^2_R \vw\right\|_{\dot{H}^{-1}}.
\end{equation}
Using the embedding $L^{\frac{6}{5}}(\R)\subset \dot{H}^{-1}(\R)$, by the H\"older inequalities (with $\frac{5}{6}=\frac{2}{3}+\frac{1}{6}$) and by the Sobolev inequalities we can write
$$\left\|\theta^2_R \vw\right\|_{\dot{H}^{-1}}\leq\left\|\theta^2_R \vw\right\|_{L^{\frac{6}{5}}}\leq \|\theta^2_R\|_{L^{\frac{3}{2}}}\|\vw\|_{L^6}\leq C_R\|\vw\|_{\dot{H}^1},$$
from which we deduce
$$\left\|\dfrac{1}{[\epsilon\Delta^2+(-\Delta) ]}\kappa(\theta^2_R \vw)\right\|_{E}\leq C_{R,\epsilon}\kappa\|\vw\|_E.$$
%%%%%%%%%%%%%%%%%%%%%%%%%%%%%%%%%%%%%%%%%%%%%%%%%%%
\item The term (7) of (\ref{Estimation_Continuite01}) can be handle just as the quantity (2) (see (\ref{Estimation_Continuite2})) and we have
$$\left\|\dfrac{1}{[\epsilon \Delta^2+(-\Delta) ]} \theta^2_R (\vn \wedge \vu)\right\|_{E}\leq C_{R, \epsilon}\|\vu\|_{E}.$$
%%%%%%%%%%%%%%%%%%%%%%%%%%%%%%%%%%%%%%%%%%%%%%%%%%%
\item Finally, the quantity (8) of (\ref{Estimation_Continuite01}) is treated in the same manner as the term (3) and we obtain
$$\left\|\dfrac{1}{[\epsilon \Delta^2+(-\Delta) ]} \vg\right\|_{E}\leq  C\|\vg\|_{\dot{H}^{-1}}+\frac{C}{\sqrt{\epsilon}}\|\vg\|_{\dot{H}^{-2}}<+\infty,$$
by the hypothesis over the external force $\vg$.\\[5mm]
\end{itemize}
With all the previous estimates for the terms (1)-(8) of (\ref{Estimation_Continuite01}), and due to the structure of the application $T_{R,\epsilon}$ (which has linear elements and bilinear ones) we finally deduce the continuity of this application and the proof of Proposition \ref{Propo_Continuity} is complete. \hfill $\blacksquare$\\

%%%%%%%%%%%%%%%%%%%%%%%%%%%%%%%%%%%%%%%%%%%%%%%%%%%
\begin{Proposition}[Compactness]\label{Propo_Compactness}
Under the general setting of the Theorem \ref{Theorem_Existence}, the application $T_{R,\epsilon}$ defined in \eqref{eq_epsilon} is compact.
\end{Proposition}
%%%%%%%%%%%%%%%%%%%%%%%%%%%%%%%%%%%%%%%%%%%%%%%%%%%
{\bf Proof.} To prove the compactness of the application $T_{R,\epsilon} $, we consider from now on a bounded sequence $( \vU_n )_{n\in \mathbb{N}}$ in $E_\sigma \times E$ and we shall prove that there exist a subsequence $T_{R,\epsilon}(\vU_{n_k})_{k\in \mathbb{N}}$ which converges strongly in $E_\sigma\times E$. To this end, we easily remark that the sequence $(\theta_R \vU_n)_{n\in \mathbb{N}}$ is also bounded in $E_\sigma \times E$ (since these spaces are based in $\dot{H}^1(\R) $ and $\dot{H}^{2}(\R)$). As $R>0$ is fixed, we can therefore assume that $supp(\theta_R \vU_n)\subset B(0,4R):=B_R$. Hence, by the Rellich-Kondrashov lemma, there exists a subsequence $(\theta_R \vU_{n_k})_{k\in \mathbb{N}}$ that converges strongly in $L^p_{loc}(\R)$ for $1\le p <6$. Therefore, in order to deduce the compactness of the operator $T_{R,\epsilon} (\vU) $, we shall bound each term of (\ref{Estimation_Continuite01}) by a suitable $L^p_{loc}$ norm. In this sense, some of the computations below are related to the ones performed above, but for the sake of completeness we give here the details.\\

Remark that since the terms (3) and (8) of (\ref{Estimation_Continuite01}) are related to the external forces, we do not need to study them here as they are not linked to the variables $\vu$ or $\vw$.
\begin{itemize}
%%%%%%%%%%%%%%%%%%%%%%%%%%%%%%%%%%%%%%%%%%%%%%%%%%%
\item For the first term of (\ref{Estimation_Continuite01}), following the same ideas which led us to (\ref{Estimate_Term1}) we obtain the estimate 
\begin{equation}\label{Estimation_Compacite}
\left\|\dfrac{1}{[\epsilon \Delta^2+(-\Delta) ]} \mathbb{P}\left[\left(\left(\theta_R \vu\right) \cdot \vn\right)\left(\theta_R \vu\right)\right]\right\|_{E_\sigma}\leq \left(C+\frac{C}{\sqrt{\epsilon}}\right)\left\|\left(\left(\theta_R \vu\right) \cdot \vn\right)\left(\theta_R \vu\right)\right\|_{\dot{H}^{-1}}.
\end{equation}
We use now the identity (using the divergence free condition for $\vu$)
$$\left(\left(\theta_R \vu\right) \cdot \vn\right)\left(\theta_R \vu\right)=div(\theta_R^2(\vu\otimes \vu))+\theta_R(\vu\cdot\vn \theta_R)\vu,$$
to write
\begin{eqnarray*}
\left\|\left(\left(\theta_R \vu\right) \cdot \vn\right)\left(\theta_R \vu\right)\right\|_{\dot{H}^{-1}}&\leq &\|div(\theta_R^2(\vu\otimes \vu))\|_{\dot{H}^{-1}}+ \|\theta_R(\vu\cdot\vn \theta_R)\vu\|_{\dot{H}^{-1}}\\
&\leq & C\|\theta_R^2(\vu\otimes \vu)\|_{L^{2}}+ \|\theta_R(\vu\cdot\vn \theta_R)\vu\|_{\dot{H}^{-1}}.
\end{eqnarray*}
By the Hölder inequalities with $\frac{1}{2} =\frac{1}{10}+\frac{1}{5}+ \frac{1}{5}$ for the first term above and by the space inclusion $L^{\frac{6}{5}}(\R)\subset \dot{H}^{-1}(\R)$ and the H\"older inequality with $\frac{5}{6} = \frac{1}{4}+\frac{1}{4}+\frac{1}{3}$ for the second quantity above, we have
\begin{eqnarray*}
&\leq & C\|\theta_R^2\|_{L^{10}}\|\vu\|_{L^{5}(B_R)}\|\vu\|_{L^{5}(B_R)}+ \|\theta_R(\vu\cdot\vn \theta_R)\vu\|_{L^{\frac{6}{5}}}\\
&\leq &C\|\theta_R^2\|_{L^{10}}\|\vu\|_{L^{5}(B_R)}\|\vu\|_{L^{5}(B_R)}+ \|\theta_R\|_{L^{3}}\|\vu\|_{L^4(B_R)}\|\vn \theta_R\|_{L^\infty}\|\vu\|_{L^{4}(B_R)},
\end{eqnarray*}
and we thus obtain the wished control with $L^p_{loc}$ norms with $1\leq p<6$.
%%%%%%%%%%%%%%%%%%%%%%%%%%%%%%%%%%%%%%%%%%%%%%%%%%%
\item For the term (2) of (\ref{Estimation_Continuite01}) we have, from (\ref{Estimate_Term2})
$$\left\|\dfrac{1}{[\epsilon \Delta^2+(-\Delta) ]} \theta^2_R (\vn \wedge \vw)\right\|_{E_\sigma}\leq C\|\theta^2_R \vw\|_{L^2}+C\| \vw \wedge \vn (\theta^2_R)\|_{L^{\frac{6}{5}}}+\frac{C}{\sqrt{\epsilon}}\|\theta^2_R \vw\|_{L^2}+\frac{C}{\sqrt{\epsilon}}\| \vw \wedge \vn (\theta^2_R)\|_{L^{\frac{6}{5}}},$$
from which we easily deduce the estimate
\begin{eqnarray}
&\leq &C\|\theta^2_R\|_{L^\infty} \|\vw\|_{L^2(B_R)}+C\| \vn (\theta^2_R)\|_{L^\infty}\| \vw\|_{L^{\frac{6}{5}}(B_R)}\notag\\
&&+\frac{C}{\sqrt{\epsilon}}\|\theta^2_R\|_{L^\infty}\|\vw\|_{L^2(B_R)}+\frac{C}{\sqrt{\epsilon}}\|\vn (\theta^2_R)\|_{L^\infty}\|\vw\|_{L^{\frac{6}{5}}(B_R)},\label{Estimation_Compacite_Terme2}
\end{eqnarray}
which is the desired inequality in terms of $L^p_{loc}$ norms with $1\leq p<6$.
%%%%%%%%%%%%%%%%%%%%%%%%%%%%%%%%%%%%%%%%%%%%%%%%%%%
\item For the term (4) of (\ref{Estimation_Continuite01}) we write
\begin{eqnarray*}
\left\|\dfrac{1}{[\epsilon \Delta^2+(-\Delta) ]}\vn (\theta_R  div(\vw))\right\|_{E}&=&\left\|\dfrac{(-\Delta)^{\frac{3}{2}}}{[\epsilon \Delta^2+(-\Delta) ]}\frac{1}{(-\Delta)^{\frac{3}{2}}}\vn (\theta_R  div(\vw))\right\|_{\dot{H}^1}\\
&&+\sqrt{\epsilon}\left\|\dfrac{\Delta^2}{[\epsilon \Delta^2+(-\Delta) ]}\frac{1}{\Delta^2}\vn (\theta_R  div(\vw))\right\|_{\dot{H}^2},
\end{eqnarray*}
and thus, by the Lemma \ref{Lemme_EstimationFourier} we obtain
\begin{eqnarray*}
&\leq &\frac{C}{\sqrt{\epsilon}}\left\|\frac{1}{(-\Delta)^{\frac{3}{2}}}\vn (\theta_R  div(\vw))\right\|_{\dot{H}^1}+\frac{C}{\sqrt{\epsilon}}\left\|\frac{1}{\Delta^2}\vn (\theta_R  div(\vw))\right\|_{\dot{H}^2}\\
&\leq &\frac{C}{\sqrt{\epsilon}}\left\|\theta_R  div(\vw)\right\|_{\dot{H}^{-1}}+\frac{C}{\sqrt{\epsilon}}\left\|\theta_R  div(\vw)\right\|_{\dot{H}^{-1}}\leq \frac{C}{\sqrt{\epsilon}}\left\|\theta_R  div(\vw)\right\|_{\dot{H}^{-1}},
\end{eqnarray*}
where in the last estimate we used the properties of the homogeneous Sobolev spaces. We use now the identity $\theta_R div(\vw)=div(\theta_R \vw)-\vn \theta_R \cdot \vw$ and we have (by the embedding $L^{\frac{6}{5}}(\R)\subset \dot{H}^{-1}(\R)$) 
\begin{eqnarray*}
\left\|\theta_R  div(\vw)\right\|_{\dot{H}^{-1}}&\leq &\left\|div(\theta_R \vw)\right\|_{\dot{H}^{-1}}+\left\|\vn \theta_R \cdot \vw\right\|_{\dot{H}^{-1}}\leq \left\|\theta_R \vw\right\|_{L^2}+\left\|\vn \theta_R \cdot \vw\right\|_{L^{\frac{6}{5}}}\\
&\leq& \|\theta_R\|_{L^\infty}\|\vw\|_{L^2(B_R)}+\|\vn \theta_R \|_{L^\infty} \|\vw\|_{L^{\frac{6}{5}}(B_R)},
\end{eqnarray*}
which is the wished control in terms of $L^p_{loc}$ norms with $1\leq p<6$.
%%%%%%%%%%%%%%%%%%%%%%%%%%%%%%%%%%%%%%%%%%%%%%%%%%%
\item The quantity (5) of (\ref{Estimation_Continuite01}) is treated as follows: by the same arguments that lead us to the estimate (\ref{Estimation_Compacite}) we can write
$$\left\|\dfrac{1}{[\epsilon \Delta^2+(-\Delta) ]}\left(\left(\theta_R \vu\right)\cdot \vn\right)\left(\theta_R \vw\right)\right\|_{E}\leq \left(C+\frac{C}{\sqrt{\epsilon}}\right)\left\|\left(\left(\theta_R \vu\right) \cdot \vn\right)\left(\theta_R \vw\right)\right\|_{\dot{H}^{-1}}.$$
By the divergence free condition of $\vu$ (recall that $\vw$ is \emph{not} divergence free), we still have the identity
$$\left(\left(\theta_R \vu\right) \cdot \vn\right)\left(\theta_R \vw\right)=div(\theta_R^2(\vw\otimes \vu))+\theta_R(\vu\cdot\vn \theta_R)\vw,$$
and we obtain
\begin{eqnarray*}
\left\|\left(\left(\theta_R \vu\right) \cdot \vn\right)\left(\theta_R \vw\right)\right\|_{\dot{H}^{-1}}&\leq &\left\|div(\theta_R^2(\vw\otimes \vu))\right\|_{\dot{H}^{-1}}+\left\|\theta_R(\vu\cdot\vn \theta_R)\vw\right\|_{\dot{H}^{-1}}\\
&\leq & C\left\|\theta_R^2(\vw\otimes \vu)\right\|_{L^2}+C\left\|\theta_R(\vu\cdot\vn \theta_R)\vw\right\|_{L^{\frac{6}{5}}},
\end{eqnarray*}
where we used the properties of the homogeneous Sobolev spaces as well as the embedding $L^{\frac{6}{5}}(\R)\subset \dot{H}^{-1}(\R)$. Now, by the Hölder inequalities with $\frac{1}{2} =\frac{1}{10}+\frac{1}{5}+ \frac{1}{5}$ and $\frac{5}{6} = \frac{1}{4}+\frac{1}{4}+\frac{1}{3}$ we obtain
$$\leq  C\|\theta_R^2\|_{L^{10}}\|\vw\|_{L^5(B_R)}\|\vu\|_{L^5(B_R)}+C\|\theta_R\|_{L^3}\|\vu\|_{L^4(B_R)}\|\vn \theta_R\|_{L^\infty}\|\vw\|_{L^{4}(B_R)},$$
and we obtain the needed estimates in terms of $L^p_{loc}$ norms with $1\leq p<6$.
%%%%%%%%%%%%%%%%%%%%%%%%%%%%%%%%%%%%%%%%%%%%%%%%%%%
\item For the term (6) of (\ref{Estimation_Continuite01}), by the estimate (\ref{Estimation_Continuite_Term6}) we easily obtain
\begin{eqnarray*}
\left\|\dfrac{1}{[\epsilon\Delta^2+(-\Delta) ]}\kappa(\theta^2_R \vw)\right\|_{E}&\leq &\kappa\left(C+\frac{C}{\sqrt{\epsilon}}\right)\left\|\theta^2_R \vw\right\|_{\dot{H}^{-1}}\leq \kappa\left(C+\frac{C}{\sqrt{\epsilon}}\right)\left\|\theta^2_R \vw\right\|_{L^{\frac{6}{5}}}\\
&\leq &\kappa\left(C+\frac{C}{\sqrt{\epsilon}}\right)\|\theta^2_R\|_{L^\infty}\| \vw\|_{L^{\frac{6}{5}}(B_R)},\end{eqnarray*}
and we have then the desired control. 
%%%%%%%%%%%%%%%%%%%%%%%%%%%%%%%%%%%%%%%%%%%%%%%%%%%
\item For the quantity (7) of (\ref{Estimation_Continuite01}), it is enough to proceed as the second term of (\ref{Estimation_Continuite01}) and to apply the same arguments used to obtain the estimate (\ref{Estimation_Compacite_Terme2}).
%%%%%%%%%%%%%%%%%%%%%%%%%%%%%%%%%%%%%%%%%%%%%%%%%%%
\end{itemize}
We have obtained suitable controls for each one of the terms of (\ref{Estimation_Continuite01}) in terms of $L^p_{loc}$ norms with $1\leq p<6$, from which we obtain the compactness of the application $T_{R,\epsilon}$. The proof of Proposition \ref{Propo_Compactness} is now complete. \hfill $\blacksquare$\\

We have proven the continuity and the compactness of the application $T_{R,\epsilon}$ and in order to apply the Theorem \ref{Teo_SchauderFixPoint}, we need now to establish some a priori estimates. 
%%%%%%%%%%%%%%%%%%%%%%%%%%%%%%%%%%%%%%%%%%%%%%%%%%%
\begin{Proposition}[A priori estimates]\label{Propo_AprioriEstimates}
Under the general setting of the Theorem \ref{Theorem_Existence}, there exists a finite $M$ such that, for every $\lambda \in[0,1]$, if $\vU = \lambda T_{R,\epsilon}(\vU)$ we then have $\|\vU\|_{E_\sigma\times E}\leq M$.
\end{Proposition}
%%%%%%%%%%%%%%%%%%%%%%%%%%%%%%%%%%%%%%%%%%%%%%%%%%%
{\bf Proof.} Since we have $\vU = \lambda T_{R,\epsilon}(\vU)$ for all $\lambda \in[0,1]$, using the expressions (\ref{FixedPoint})-(\ref{eq_epsilon}) we can write

\begin{eqnarray*}
\vU&=&\lambda T_{R,\epsilon} (\vU)\\
\left(\begin{array}{l} 
\vu\\[5mm]
\vw
\end{array}\right)& = &\dfrac{\lambda}{[\epsilon \Delta^2  + (-\Delta) ]} 
\left(\begin{array}{c}  
-\mathbb{P}\left[\left(\left(\theta_R \vu\right) \cdot \vn\right)\left(\theta_R \vu\right)\right]
+ \frac{1}{2}\theta^2_R (\vn \wedge \vw)+\vf\\ [5mm]
\vn (\theta_R  div(\vw))-\left(\left(\theta_R \vu\right) \cdot \vn\right)\left(\theta_R \vw\right)-\kappa(\theta^2_R \vw)+\frac{1}{2}\theta^2_R (\vn \wedge \vu)+  \vec{g}
\end{array}\right),
\end{eqnarray*}
from which we deduce the system
\begin{eqnarray*} 
[\epsilon \Delta^2  + (-\Delta) ]\vu & =&-\lambda\left(\mathbb{P}\left[\left(\left(\theta_R \vu\right) \cdot \vn\right)\left(\theta_R \vu\right)\right]+ \frac{1}{2}\theta^2_R (\vn \wedge \vw)+\vf\right)\\ [3mm]
[\epsilon \Delta^2  + (-\Delta) ]\vw&=&\lambda\left(\vn (\theta_R  div(\vw))-\left(\left(\theta_R \vu\right) \cdot \vn\right)\left(\theta_R \vw\right)-\kappa(\theta^2_R \vw)+\frac{1}{2}\theta^2_R (\vn \wedge \vu)+  \vec{g}\right).
\end{eqnarray*}
Now, if we multiply the first equation above by $\vu$, the second by $\vw$ and we integrate over $\R$ we can write (after an integration by parts in the left-hand side):
\begin{eqnarray*} 
\epsilon\int_{\R}|\Delta \vu|^2dx+\int_{\R}|\vn\otimes \vu|^2dx& =&-\lambda\int_{\R}\mathbb{P}\left[\left(\left(\theta_R \vu\right) \cdot \vn\right)\left(\theta_R \vu\right)\right]\cdot \vu\,dx+ \frac{\lambda}{2}\int_{\R}\theta^2_R (\vn \wedge \vw)\cdot \vu\, dx\\
&&+\lambda \int_{\R}\vf\cdot \vu\, dx\\ [3mm]
\epsilon\int_{\R}|\Delta \vw|^2dx+\int_{\R}|\vn\otimes \vw|^2dx&=&\lambda\int_{\R}\vn (\theta_R  div(\vw))\cdot \vw \, dx-\lambda \int_{\R}\left(\left(\theta_R \vu\right) \cdot \vn\right)\left(\theta_R \vw\right)\cdot \vw\, dx\\
&&-\lambda \int_{\R}\kappa(\theta^2_R \vw)\cdot \vw\, dx+\frac{\lambda}{2} \int_{\R}\theta^2_R (\vn \wedge \vu)\cdot \vw\, dx+ \lambda\int_{\R}\vec{g}\cdot\vw\, dx.
\end{eqnarray*}
We note now that, since $div(\vu)=0$ and by the properties of the Leray projector $\mathbb{P}$, we have
$$\int_{\R}\mathbb{P}[(\theta_R\vu\cdot \grad)(\theta_R\vu)] \cdot \vu\, dx=\int_{\R}(\theta_R\vu\cdot \grad)(\theta_R\vu) \cdot \vu\, dx=0\quad \mbox{and}\quad \int_{\R}(\theta_R\vu\cdot \grad)(\theta_R\vw) \cdot \vw \,dx=0,$$
thus, by adding the two equations we obtain
\begin{eqnarray*} 
&&\epsilon\left(\int_{\R}|\Delta \vu|^2dx+\int_{\R}|\Delta \vw|^2dx\right)+\int_{\R}|\vn\otimes \vu|^2dx+\int_{\R}|\vn\otimes \vw|^2dx\\
&&=\frac{\lambda}{2}\left(\int_{\R}\theta^2_R (\vn \wedge \vw)\cdot \vu\, dx+ \int_{\R}\theta^2_R (\vn \wedge \vu)\cdot \vw\, dx\right)+\lambda\int_{\R}\vn (\theta_R  div(\vw))\cdot \vw \, dx-\lambda \int_{\R}\kappa(\theta^2_R \vw)\cdot \vw\, dx\\
&&+\lambda \int_{\R}\vf\cdot \vu\, dx+ \lambda\int_{\R}\vec{g}\cdot\vw\, dx.
\end{eqnarray*}
Now, by an integration by parts we observe that we have the identity 
\begin{eqnarray*}
\int_{\R}\theta^2_R(\rot \vw)\cdot \vu\, dx&=&\int_{\R}\vw \cdot \rot(\theta_R^2 \vu) dx\\
&=&\int_{\R}\theta^2_R (\rot \vu)\cdot \vw dx+\int_{\R}\vn(\theta^2_R)\cdot (\vw\wedge \vu) dx,
\end{eqnarray*}
and we can write
\begin{eqnarray*} 
\epsilon\left(\| \vu\|_{\dot{H}^2}^2+ \|\vw\|_{\dot{H}^2}^2\right)+\|\vu\|_{\dot{H}^1}^2+\|\vw\|_{\dot{H}^1}^2&=&\lambda\int_{\R}\theta^2_R (\rot \vu)\cdot \vw \,dx+\frac{\lambda}{2}\int_{\R}\theta_R\vn(\theta_R)\cdot (\vw\wedge \vu) dx \\&&+\lambda\int_{\R}\vn (\theta_R  div(\vw))\cdot \vw \, dx
-\lambda \kappa\int_{\R}\theta^2_R |\vw|^2 dx\\
&&+\lambda \int_{\R}\vf\cdot \vu\, dx+ \lambda\int_{\R}\vec{g}\cdot\vw\, dx.
\end{eqnarray*}
Noting that $\displaystyle{\int_{\R}\vn (\theta_R  div(\vw))\cdot \vw \, dx=-\int_{\R}\theta_R  |div(\vw)|^2dx}$, we have
\begin{eqnarray*} 
\epsilon\left(\| \vu\|_{\dot{H}^2}^2+ \|\vw\|_{\dot{H}^2}^2\right)+\|\vu\|_{\dot{H}^1}^2+\|\vw\|_{\dot{H}^1}^2&\leq& \lambda \int_{\R}\theta^2_R (\rot \vu)\cdot \vw \,dx+C\lambda\int_{\R}|\vn(\theta_R)||\theta_R\vw||\vu| dx \\
&&-\lambda\int_{\R}\theta_R  |div(\vw)|^2dx-\lambda \kappa\int_{\R}\theta^2_R |\vw|^2 dx\\
&&+\lambda \int_{\R}\vf\cdot \vu\, dx+ \lambda\int_{\R}\vec{g}\cdot\vw\, dx,
\end{eqnarray*}
which can be rewritten as
\begin{eqnarray*} 
\epsilon\left(\| \vu\|_{\dot{H}^2}^2+ \|\vw\|_{\dot{H}^2}^2\right)+\|\vu\|_{\dot{H}^1}^2+\|\vw\|_{\dot{H}^1}^2+\lambda\int_{\R}\theta_R  |div(\vw)|^2dx\leq \lambda\int_{\R}|\theta^2_R| |\rot \vu| | \vw|dx\\
+C\lambda\int_{\R}|\vn(\theta_R)||\theta_R\vw||\vu| dx -\lambda \kappa\int_{\R}\theta^2_R |\vw|^2 dx+\lambda \int_{\R}\vf\cdot \vu\ dx+ \lambda\int_{\R}\vec{g}\cdot\vw\, dx,
\end{eqnarray*}
but since $\displaystyle{\int_{\R}\theta_R  |div(\vw)|^2dx}\geq 0$, we have
\begin{eqnarray*} 
\epsilon\left(\| \vu\|_{\dot{H}^2}^2+ \|\vw\|_{\dot{H}^2}^2\right)+\|\vu\|_{\dot{H}^1}^2+\|\vw\|_{\dot{H}^1}^2&\leq & \lambda\|\theta_R\|_{L^\infty}\int_{\R}|\theta_R| |\rot \vu| | \vw| dx+C\lambda\int_{\R}|\vn(\theta_R)||\theta_R\vw||\vu| dx \\
&&-\lambda \kappa\int_{\R}\theta^2_R |\vw|^2 dx+\lambda \int_{\R}\vf\cdot \vu dx+ \lambda\int_{\R}\vec{g}\cdot\vw\, dx.
\end{eqnarray*}
At this point, recalling that $\|\theta_R\|_{L^\infty}\le 1$, by the Cauchy-Schwarz inequality, the H\"older inequalities with $1=\frac{1}{3}+\frac{1}{2}+\frac{1}{6}$ as well as the  $\dot{H}^{-1}-\dot{H}^{1}$ duality we obtain the estimate
\begin{eqnarray*} 
\epsilon\left(\| \vu\|_{\dot{H}^2}^2+ \|\vw\|_{\dot{H}^2}^2\right)+\|\vu\|_{\dot{H}^1}^2+\|\vw\|_{\dot{H}^1}^2&\leq& \lambda\|\theta_R\vw\|_{L^2}\|\rot\vu\|_{L^2}+C\lambda\|\vn(\theta_R)\|_{L^3}\|\theta_R\vw\|_{L^2}\|\vu\|_{L^6}\\
&&-\lambda \kappa\|\theta_R \vw\|^2_{L^2}+\lambda\|\vf\|_{\dot{H}^{-1}}\|\vu\|_{\dot{H}^{1}}+ \lambda\|\vg\|_{\dot{H}^{-1}}\|\vw\|_{\dot{H}^{1}}.
\end{eqnarray*}
Moreover, since $\|\rot \vu\|_{L^2}\le 2\|\vu\|_{\dot H^1}$,  by the Young inequality for the sum, we have for some small technical parameter $0<\varepsilon<\frac{1}{3}$:
\begin{eqnarray*}
&\epsilon\left(\| \vu\|_{\dot{H}^2}^2+ \|\vw\|_{\dot{H}^2}^2\right)+\|\vu\|_{\dot{H}^1}^2+\|\vw\|_{\dot{H}^1}^2\leq& C(\varepsilon)\lambda \|\theta_R\vw\|_{L^2}^2+ \varepsilon\lambda\|\vu\|_{\dot H^1}^2+ C\lambda\|\vn(\theta_R)\|_{L^3}\|\theta_R\vw\|_{L^2}\|\vu\|_{\dot{H}^1}\\
&&-\lambda \kappa\|\theta_R \vw\|^2_{L^2}+\lambda\|\vf\|_{\dot{H}^{-1}}\|\vu\|_{\dot{H}^{1}}+ \lambda\|\vg\|_{\dot{H}^{-1}}\|\vw\|_{\dot{H}^{1}}.
\end{eqnarray*}
where we used the Sobolev embedding $L^{6}(\R)\subset \dot{H}^1(\R)$ in the estimate above.
We observe now that by homogeneity we have $\|\vn(\theta_R)\|_{L^3}=\|\vn\theta\|_{L^3}$, and using the Young inequalities for the sum we obtain
\begin{eqnarray*} 
\epsilon\left(\| \vu\|_{\dot{H}^2}^2+ \|\vw\|_{\dot{H}^2}^2\right)+\|\vu\|_{\dot{H}^1}^2+\|\vw\|_{\dot{H}^1}^2&\leq&C(\varepsilon)\lambda \|\theta_R\vw\|_{L^2}^2+ \varepsilon\lambda\|\vu\|_{\dot H^1}^2+C(\varepsilon)\lambda\|\vn\theta\|_{L^3}^2\|\theta_R\vw\|_{L^2}^2\\
&&+\varepsilon\lambda\|\vu\|_{\dot{H}^1}^2-\lambda \kappa\|\theta_R \vw\|^2_{L^2}+C(\varepsilon)\lambda\|\vf\|_{\dot{H}^{-1}}^2+ \varepsilon \lambda\|\vu\|_{\dot{H}^{1}}^2\\
&&+ C\lambda\|\vg\|_{\dot{H}^{-1}}^2+\frac{\lambda}{4}\|\vw\|_{\dot{H}^{1}}^4
\end{eqnarray*}
$$\leq \lambda\left(C(\varepsilon)+C(\varepsilon)\|\vn\theta\|_{L^3}^2-\kappa\right)\|\theta_R\vw\|_{L^2}^2+3\varepsilon \lambda\|\vu\|_{\dot{H}^1}^2+C(\varepsilon)\lambda\|\vf\|_{\dot{H}^{-1}}^2+ C\lambda\|\vg\|_{\dot{H}^{-1}}^2+\frac{\lambda}{4}\|\vw\|_{\dot{H}^{1}}^4,$$
We thus obtain the inequality
\begin{eqnarray} 
\epsilon\left(\| \vu\|_{\dot{H}^2}^2+ \|\vw\|_{\dot{H}^2}^2\right)+(1-3\varepsilon)\|\vu\|_{\dot{H}^1}^2+(1-\frac{\lambda}{4})\|\vw\|_{\dot{H}^1}^2+\lambda\left(\kappa-C(\varepsilon)-C(\varepsilon)\|\vn\theta\|_{L^3}^2\right)\|\theta_R\vw\|_{L^2}^2\notag\\
\leq C(\varepsilon)\lambda\|\vf\|_{\dot{H}^{-1}}^2+ C\lambda\|\vg\|_{\dot{H}^{-1}}^2,\label{Estimation_Uniforme_NormeL2}
\end{eqnarray}
but since $0<\varepsilon<\frac{1}{3}$, $0\leq \lambda\leq 1$, and $\kappa-C(\varepsilon)-C(\varepsilon)\|\vn\theta\|_{L^3}^2\geq 0$ as $\kappa \gg 1$   we have
\begin{equation}\label{Estimation_Uniforme}
\epsilon\| \vu\|_{\dot{H}^2}^2+ \epsilon\|\vw\|_{\dot{H}^2}^2+\|\vu\|_{\dot{H}^1}^2+\|\vw\|_{\dot{H}^1}^2\leq C\left(\|\vf\|_{\dot{H}^{-1}}^2+\|\vg\|_{\dot{H}^{-1}}^2\right),
\end{equation}
from which we deduce (since $\vf, \vg\in \dot{H}^{-1}(\R)$ by hypothesis) the estimate 
$$\|\vU\|^2_{E_\sigma\times E}\leq M<+\infty.$$
The proof of Proposition \ref{Propo_AprioriEstimates} is thus finished. \hfill$\blacksquare$\\

With the help of the Propositions \ref{Propo_Continuity}, \ref{Propo_Compactness} and \ref{Propo_AprioriEstimates}, we have verified all the hypotheses needed in order to apply the Theorem \ref{Teo_SchauderFixPoint}: if we apply this result to the system (\ref{FixedPoint}) we thus obtain the existence of a solution $(\vu_{R, \epsilon}, \vw_{R, \epsilon})$ (depending on $R$ and $\epsilon$) which satisfies the system (\ref{Probleme_Approx}) as well as the estimate (\ref{Estimation_Uniforme}).\\

It is important to remark here that the estimates above are an important source of information: indeed, we have the following result:
\begin{Corollaire}[Energy a priori estimate for $\vw_{R, \epsilon}$]\label{Coro_L2Norm}
In the general framework of the Theorem \ref{Theorem_Existence}, and for the variable  $\vw_{R, \epsilon}$ obtained above, we have for all $R>1$ and for all $0<\epsilon\leq 1$, the uniform control 
$$\|\theta_R\vw_{R, \epsilon}\|_{L^2}^2\leq C(\|\vf\|_{\dot{H}^{-1}}^2+\|\vg\|_{\dot{H}^{-1}}^2).$$
\end{Corollaire}
This is a direct consequence of the inequality (\ref{Estimation_Uniforme_NormeL2}) and we will see later on how to exploit this uniform estimate.
%%%%%%%%%%%%%%%%%%%%%%%%%%%%%%%%%%%%%%%%%%%%%%%%%%%
\paragraph{End of the proof of Theorem \ref{Theorem_Existence}.} We shall now recover the initial problem (\ref{equationwithoutpressure}) by making $R\to +\infty$ and $\epsilon\to 0$. For this, we will first fix $\epsilon>0$ and we will take the limit $R\to +\infty$: indeed, since we have at our disposal the uniform in $R$ estimate (\ref{Estimation_Uniforme}), we can extract a subsequence $R_k\to +\infty$ such that $(\vu_{R_k, \epsilon}, \vw_{R_k, \epsilon})$ converges weakly-$*$ in $\dot{H}^1(\R)$ to some limit $(\vu_{\epsilon}, \vw_{\epsilon})$ (by the Banach-Aloaglu theorem). Moreover, by the Rellich-Kondrachov lemma, we also have the strong convergence of $(\vu_{R_k, \epsilon}, \vw_{R_k, \epsilon})$ to some limit $(\vu_{\epsilon}, \vw_{\epsilon})$ in the space $L^p_{loc}(\R)$ with $1\leq p<6$. These arguments  allows us to obtain a weak convergence (in $\mathcal{D}'$) of the terms $\mathbb{P}\big[\big((\theta_R\vu_{R_k, \epsilon}) \cdot \vn\big) (\theta_R\vu_{R_k, \epsilon})\big]$ and $\big((\theta_R\vu_{R_k, \epsilon}) \cdot \vn\big) (\theta_R\vw_{R_k, \epsilon})$ to $\mathbb{P}\big[\big(\vu_{\epsilon}\cdot \vn\big) \vu_{\epsilon}\big]$ and  $\big(\vu_{\epsilon} \cdot \vn\big) \vw_{\epsilon}$ when $R\to +\infty$, respectively. We thus obtain a vector $(\vu_{\epsilon}, \vw_{\epsilon})$ which is a solution of the problem
$$
\begin{cases}
0=-\epsilon \Delta^2 \vu_{\epsilon} + \Delta \vu_{\epsilon}-\mathbb{P}\big[\big(\vu_{\epsilon} \cdot \vn\big) \vu_{\epsilon}\big]+\frac{1}{2} \vn \wedge \vw_{\epsilon}+\vf,
\\[5mm]
0=-\epsilon \Delta^2 \vw_{\epsilon} + \Delta \vw_{\epsilon}+\vn div(\vw_{\epsilon})-\big(\vu_{\epsilon} \cdot \vn\big) \vw_{\epsilon}-\kappa\vw_{\epsilon}+\frac{1}{2}\vn \wedge \vu_{\epsilon}+\vg.
\end{cases}
$$
By the same arguments as before, since we still have the uniform (in $\epsilon$) inequality $\|\vu_{\epsilon}\|_{\dot{H}^1}^2+\|\vw_{\epsilon}\|_{\dot{H}^1}^2\leq C(\|\vf\|_{\dot{H}^{-1}}^2+\|\vg\|_{\dot{H}^{-1}}^2),$ given in the estimate (\ref{Estimation_Uniforme}), there exists a subsequence $\epsilon_k\to 0$ such that $(\vu_{\epsilon_k}, \vw_{\epsilon_k})$ convergence weakly to a limit $(\vu, \vw)$ in the space $\dot{H}^1(\R)$. Thus, again by the Rellich-Kondrachov lemma, we obtain the strong convergence of $(\vu_{\epsilon_k}, \vw_{\epsilon_k})$ to $(\vu, \vw)$ in $L^p_{loc}(\R)$ with $1\leq p<6$ and from these facts we deduce the weak convergence (in $\mathcal{D}'$) of the quantities $\mathbb{P}\big[\big(\vu_{\epsilon} \cdot \vn\big) \vu_{\epsilon}\big]$ and $\big(\vu_{\epsilon} \cdot \vn\big) \vw_{\epsilon}$ to  $\mathbb{P}\big[\big(\vu \cdot \vn\big) \vu\big]$ and $\big(\vu \cdot \vn\big) \vw$
 when $\epsilon\to 0$, respectively. We have thus obtained a solution $(\vu, \vw)$ of the initial equation (\ref{equationwithoutpressure}) such that $\vu, \vw\in \dot{H}^{1}(\R)$.
 
Let us remark now that for the variable $\vw$ we have some additional information. Indeed, by the Corollary \ref{Coro_L2Norm} above we have a $L^2$-norm uniform estimate (in $R$ and $\epsilon$) of the terms $\vw_{R, \epsilon}$ from which we obtain that, up to a subsequence, we have $\vw_{R, \epsilon}\longrightarrow \vw$ in $L^2$. This fact allows us to deduce that $\vw\in L^2(\R)\cap  \dot{H}^{1}(\R)$ which is $\vw\in H^{1}(\R)$.

To end the proof, we need to study the pressure $p$ and from the expression (\ref{Equation_Pression_intro}) we obtain the identity 
\begin{equation}\label{Formule_Pression}
p=\frac{1}{(-\Delta)} div (div(\vu\otimes\vu)). 
\end{equation}
We can thus write
\begin{align*}
\left\|p\right\|_{\dot H^{\frac{1}{2}}}\le
\left\|\frac{1}{(-\Delta)} div div(\vu\otimes\vu)\right\|_{\dot H^{\frac{1}{2}}}\le \left\|\vu\otimes\vu\right\|_{\dot H^{\frac{1}{2}}},
\end{align*}
and using the usual product laws in Sobolev spaces\footnote{Recall that, for $0\leq s<+\infty$ and $0<\delta<\frac32$, we have
$\|fg\|_{\dot{H}^{s+\delta-\frac32}}\leq C\left(\|f\|_{\dot{H}^{\delta}}\|g\|_{\dot{H}^{s}}+\|g\|_{\dot{H}^{\delta}}\|f\|_{\dot{H}^{s}}\right).$ See \cite[Lemma 7.3]{lemarie2018navier} for a proof of this inequality.} we have $\left \|\vu\otimes\vu\right\|_{\dot H^{\frac{1}{2}}}\le \left\|\vu\right\|^2_{\dot H^{1}}$, from which we easily deduce that $p\in \dot H^{\frac{1}{2}}(\R)$. The proof of Theorem \ref{Theorem_Existence} is finished. \hfill $\blacksquare$
%%%%%%%%%%%%%%%%%%%%%%%%%%%%%%%%%%%%%%%%%%%%%%%%%%%
\section{Proof of Theorem \ref{Theorem_Regularity}}\label{Secc_Theo2}
We have just obtained solutions $(\vu, p, \vw)$ of the system (\ref{SteadyMicropolarEquations}) such that $\vu\in \dot{H}^1(\R)$ and $\vw\in H^1(\R)$. We are going to see now how it is possible to improve this information in the setting of the Theorem \ref{Theorem_Regularity}, \emph{i.e.} when $\vf,\vg$ are regular external forces or when $\vf=\vg=0$. For the sake of simplicity and without any loss of generality, we will assume from now on that $\vf=\vg=0$.\\

First note that due to the expression (\ref{Equation_Pression_intro}) we can restrict our study to the variables $\vu$ and $\vw$. Thus, by applying the Leray projector to the first equation of (\ref{SteadyMicropolarEquations}), with the divergence free property of $\vu$ and $\rot \vw$, we obtain 
\begin{equation}\label{Equation_U_Regularite}
\vu=\frac{1}{(-\Delta)}\left(-\mathbb{P}(div(\vu\otimes\vu))+\frac{1}{2}\rot \vw\right),
\end{equation}
and we can write
\begin{eqnarray*}
\|\vu\|_{\dot H^{\frac{3}{2}}}&\leq& \left\|\frac{1}{(-\Delta)}\mathbb{P}(div(\vu\otimes\vu))\right\|_{\dot H^{\frac{3}{2}}}
+\frac{1}{2}\left\|\frac{1}{(-\Delta)}\rot \vw\right\|_{\dot H^{\frac{3}{2}}}\\
&\leq & C\|\vu\otimes\vu\|_{\dot H^{\frac{1}{2}}}+C\|\vw\|_{\dot H^{\frac{1}{2}}},
\end{eqnarray*}
where we used the boundedness properties of the Leray projector as well as the properties of homogeneous Sobolev spaces. At this point we remark that, as $\vw\in L^2(\R)\cap \dot{H}^1(\R)$ we also have $\vw\in \dot{H}^{\frac{1}{2}}(\R)$ (with $\|\vw\|_{\dot{H}^{\frac{1}{2}}}\leq \|\vw\|_{L^2}^{\frac{1}{2}}\|\vw\|_{\dot{H}^1}^{\frac{1}{2}}$) and by the product laws in Sobolev space, we can write
$$\|\vu\|_{\dot H^{\frac{3}{2}}}\leq  C\|\vu\|_{\dot H^{1}}\|\vu\|_{\dot H^{1}}+C\|\vw\|_{L^2}^{\frac{1}{2}}\|\vw\|_{\dot{H}^1}^{\frac{1}{2}}<+\infty,$$
from which we deduce a first gain of regularity for $\vu$ as we have now $\vu\in \dot{H}^{\frac{3}{2}}(\R)$. Plugging this information in (\ref{Equation_U_Regularite}) we can now study a second gain of regularity for $\vu$. Indeed we have, by the same arguments as above:
\begin{eqnarray}
\|\vu\|_{\dot H^{2}}&\leq &\left\|\frac{1}{(-\Delta)}\mathbb{P}(div(\vu\otimes\vu))\right\|_{\dot H^{2}}
+\frac{1}{2}\left\|\frac{1}{(-\Delta)}\rot \vw\right\|_{\dot H^{2}}\notag\\
&\leq &C\|\vu\otimes\vu\|_{\dot H^{1}}+C\|\vw\|_{\dot H^{1}}\leq C\|\vu\|_{\dot H^{1}}\|\vu\|_{\dot H^{\frac{3}{2}}}+C\|\vw\|_{\dot H^{1}}<+\infty.\label{Formule_Dependence_ReguUW}
\end{eqnarray}
We have then $\vu\in \dot{H}^1(\R)\cap\dot{H}^2(\R)$ from which we can deduce (by working in the Fourier level) that $\vu\in L^\infty(\R)$.\\

As it is clear from the estimate (\ref{Formule_Dependence_ReguUW}), the regularity of $\vu$ is linked to the regularity of $\vw$ and in order to get more regularity for $\vu$ we will need to study the regularity of the variable $\vw$. Thus, having in mind that $\vu \in \dot{H}^1(\R)\cap\dot{H}^2(\R)\cap L^\infty(\R)$ we will now study the variable $\vw\in H^1(\R)$. However, it will be more convenient to consider first the quantity $div(\vw)$: indeed, by applying the divergence operator to the second equation of  (\ref{SteadyMicropolarEquations}) and since we have $div(\vn \wedge \vu)=0$, we obtain:
$$2\Delta div(\vw)=\kappa div(\vw)+div((\vu \cdot \vn) \vw),$$
from which we can write 
\begin{eqnarray*}
\|div(\vw)\|_{\dot{H}^{1}}&\leq &C\kappa \left\|\frac{1}{(-\Delta)}div(\vw)\right\|_{\dot{H}^1}+C\left\|\frac{1}{(-\Delta)}div((\vu \cdot \vn) \vw)\right\|_{\dot{H}^1}\\
&\leq & C\kappa \|\vw\|_{L^2}+C\|(\vu \cdot \vn) \vw\|_{L^2}\leq C\kappa \|\vw\|_{L^2}+C\|\vu\|_{L^\infty} \|\vw\|_{\dot{H}^1}<+\infty,
\end{eqnarray*}
and we obtain that $div(\vw)\in \dot{H}^{1}(\R)$. With this preliminary information over $div(\vw)$, we come back to the study of the variable $\vw$ and by the second equation of (\ref{SteadyMicropolarEquations}) we deduce the equation
$$\vw=\frac{1}{(-\Delta)}\vn div(\vw)-\kappa\frac{1}{(-\Delta)}\vw-\frac{1}{(-\Delta)}(\vu \cdot \vn) \vw+\frac{1}{2}\frac{1}{(-\Delta)} \vn \wedge \vu,$$
and we write
\begin{eqnarray*}
\|\vw\|_{\dot{H}^2}&\leq &\left\|\frac{1}{(-\Delta)}\vn div(\vw)\right\|_{\dot{H}^2}+\kappa\left\|\frac{1}{(-\Delta)}\vw\right\|_{\dot{H}^2}+\left\|\frac{1}{(-\Delta)}(\vu \cdot \vn) \vw\right\|_{\dot{H}^2}+C\left\|\frac{1}{(-\Delta)} \vn \wedge \vu\right\|_{\dot{H}^2}\\
&\leq& C\left\|div(\vw)\right\|_{\dot{H}^1}+C\kappa\|\vw\|_{L^2}+C\|(\vu \cdot \vn) \vw\|_{L^2}+C\|\vn \wedge \vu\|_{L^2}\\
&\leq & C\left\|div(\vw)\right\|_{\dot{H}^1}+C\kappa\|\vw\|_{L^2}+C\|\vu\|_{L^\infty}\|\vw\|_{\dot{H}^1}+C\|\vu\|_{\dot{H}^1}<+\infty,
\end{eqnarray*}
we thus have obtained that $\vw\in \dot{H}^2(\R)$.\\ 

Now, this gain of regularity for the variable $\vw$ will cause a gain of regularity for the variable $\vu$:  indeed, by the same arguments displayed to obtain (\ref{Formule_Dependence_ReguUW}) we can prove that $\vu\in \dot{H}^3(\R)$. But since we have $\vn\otimes \vu\in \dot{H}^1(\R)\cap \dot{H}^2(\R)$ we also have that $\vn\otimes \vu\in L^\infty(\R)$ and then, as before, we obtain that $div(\vw)\in \dot{H}^2(\R)$, from which we easily deduce that $\vw\in \dot{H}^3(\R)$. Then, by bootstrapping the same ideas as above, we will obtain that the variables $\vu$ and $\vw$ are regular and this ends the proof of Theorem \ref{Theorem_Regularity}.\hfill $\blacksquare$
%%%%%%%%%%%%%%%%%%%%%%%%%%%%%%%%%%%%%%%%%%%%%%%%%%%
\section{Proof of Theorem \ref{Theorem_Liouville_type}}\label{Secc_Theo3}
We have now $(\vu, p, \vw)$ a weak solution of the system (\ref{SteadyMicropolarEquations}) such that $\vu\in \dot{H}^{1}(\R)$, $p\in \dot{H}^{\frac{1}{2}}(\R)$ and $\vw\in H^1(\R)$. Moreover we have also, by hypothesis, the information $\vu\in L^q(\R)$ with $3\leq q\leq \frac{9}{2}$, from which we easily deduce that $p\in L^{\frac{p}{2}}(\R)$, indeed, we can write by the expression (\ref{Formule_Pression}):
\begin{equation}\label{Estimation_Pression_Lq2}
\|p\|_{L^\frac{p}{2}}=\left\|\frac{1}{(-\Delta)} div (div(\vu\otimes\vu))\right\|_{L^\frac{p}{2}}\leq C\|\vu\otimes\vu\|_{L^\frac{p}{2}}\leq C\|\vu\|_{L^p}\|\vu\|_{L^p}<+\infty.
\end{equation}
With all this information for $\vu$, $p$ and $\vw$ we will prove that if the external forces are null, then the trivial solution is unique. To do so, we consider $\phi$ a positive, smooth cut-off function such that $0\leq \phi\leq1,$ $\phi(x)=1$ if $|x|<1$ and $\phi(x)=0$ if $|x|>2$. Moreover, for $R\geq 1$ we define the function $\phi_R(\cdot)=\phi(\frac{\cdot}{R})$, note in particular that $sup(\phi_R)=B(0,R):=B_R$.\\

Multiplying the two equations of the system \eqref{SteadyMicropolarEquations} by $\phi_R^2 \vu$ and $\phi_R^2 \vw$, respectively, we obtain (recall that we are considering here $\vf=\vg=0$):
$$\begin{cases}
\Delta \vu\cdot(\phi_R^2 \vu)-[(\vu \cdot \vn) \vu]\cdot(\phi_R^2 \vu)-\vn p\cdot(\phi_R^2 \vu)+\frac{1}{2} [\vn \wedge \vw(\phi_R^2 \vu)]\cdot (\phi_R^2 \vu)=0,\\[5mm]
\Delta \vw\cdot(\phi_R^2 \vw)+[\vn div(\vw)]\cdot(\phi_R^2 \vw)-\kappa\vw\cdot(\phi_R^2 \vw)-[(\vu \cdot \vn) \vw]\cdot(\phi_R^2 \vw)+\frac{1}{2} [\vn \wedge \vu]\cdot(\phi_R^2 \vw)=0.
\end{cases}$$
Now, integrating over $\R$ and by adding these two equations we can write, after some rearrangements:
\begin{eqnarray}
&& -\int_{\R}\Delta \vu \cdot (\phi_R^2 \vu) dx-\int_{\R}\Delta \vw \cdot (\phi_R^2 \vw)dx-\int_{\R}[\grad div(\vw)]\cdot \phi_R^2 \vw dx + \kappa\int_{\R} \phi_R^2|\vw|^2dx\qquad  \notag\\
 &=& -\int_{\R}[(\vu\cdot \grad)\vu] \cdot (\phi_R^2 \vu) dx-\int_{\R}[(\vu\cdot \grad)\vw]\cdot (\phi_R^2\vw) dx - \int_{\R}\grad p\cdot (\phi_R^2 \vu) dx\notag\\
&&+\frac{1}{2}\int_{\R}\rot \vw\cdot (\phi_R^2 \vu) dx+\frac{1}{2}\int_{\R}\rot \vu \cdot (\phi_R^2 \vw) dx.\label{LocalEstimateOmega}
\end{eqnarray}
Since by the Theorem \ref{Theorem_Regularity} we have enough regularity for the functions $\vu, p$ and $\vw$, by an integration by parts, we obtain the following identities for the three first terms above
\begin{eqnarray*}
-\int_{\R}\Delta \vu \cdot(\phi_R^2 \vu) dx &= &\int_{\R} \phi_R^2| \grad\otimes \vu|^2 dx- \frac{1}{2}\int_{\R}\Delta(\phi_R^2) |\vu|^2dx\\
-\int_{\R}\Delta \vw \cdot(\phi_R^2 \vw) dx &= &\int_{\R} \phi_R^2| \grad\otimes \vw|^2 dx- \frac{1}{2}\int_{\R}\Delta(\phi_R^2) |\vw|^2dx,
\end{eqnarray*}
and
\begin{equation*}
-\int_{\R}[\grad div(\vw)]\cdot (\phi_R^2 \vw )dx=\int_{\R}  \phi_R ^2|div( \vw)|^2 dx+\int_{\R}\grad (\phi_R^2)\cdot  \vw div (\vw)dx.
\end{equation*}
Hence replacing the previous identities in \eqref{LocalEstimateOmega} we have
\begin{align*}
&\int_{\R} \phi_R^2| \grad\otimes \vu|^2dx+\int_{\R} \phi_R^2| \grad\otimes \vw|^2dx+\int_{\R} \phi_R^2| div( \vw)|^2 dx+\kappa \int_{\R} \phi_R^2|\vw|^2dx\notag\\
&=\frac{1}{2}\int_{\R}\Delta(\phi_R^2) (|\vu|^2+|\vw|^2) dx-\int_{\R}\grad (\phi_R^2)\cdot \vw div (\vw) dx-\int_{\R}[(\vu\cdot \grad)\vu]\cdot (\phi_R^2\vu) dx\notag\\
&-\int_{\R}[(\vu\cdot \grad)\vw] \cdot (\phi_R^2 \vw) dx-\int_{\R}\grad p\cdot (\phi_R^2 \vu) dx\\
&+\frac{1}{2}\int_{\R}\rot \vw\cdot( \phi_R^2 \vu) dx +\frac{1}{2}\int_{\R} \rot\vu \cdot (\phi_R^2 \vw) dx.\notag
\end{align*}
Additionally, by integration by parts we observe that we have the identity
\begin{eqnarray*}
\frac{1}{2}\int_{\R}\rot \vw\cdot( \phi_R^2\vu) dx &=&\frac{1}{2}\int_{\R} \vw\cdot\rot( \phi_R^2 \vu) dx  \\
&=&\frac{1}{2}\int_{\R} \vw\cdot[\vn(\phi_R^2)\wedge \vu] dx+\frac{1}{2}\int_{\R} \rot\vu \cdot (\phi_R^2 \vw) dx.
\end{eqnarray*}
Therefore, we can write
\begin{align}
&\int_{\R} \phi_R^2| \grad\otimes \vu|^2dx+\int_{\R} \phi_R^2| \grad\otimes \vw|^2dx+\int_{\R} \phi_R^2| div( \vw)|^2 dx+\kappa \int_{\R} \phi_R^2|\vw|^2dx\notag\\
&=\frac{1}{2}\underbrace{\int_{\R}\Delta(\phi_R^2) (|\vu|^2+|\vw|^2) dx}_{(1)}-\underbrace{\int_{\R}\grad (\phi_R^2)\cdot \vw div (\vw) dx}_{(2)}-\underbrace{\int_{\R}[(\vu\cdot \grad)\vu]\cdot (\phi_R^2\vu) dx}_{(3)}\notag\\
&-\underbrace{\int_{\R}[(\vu\cdot \grad)\vw] \cdot (\phi_R^2 \vw) dx}_{(4)}-\underbrace{ \int_{\R}\grad p\cdot (\phi_R^2 \vu) dx}_{(5)}\label{LocalEstimateBothVariables}\\
&+\underbrace{\frac{1}{2}\int_{\R} \vw\cdot[\vn(\phi_R^2)\wedge \vu] dx}_{(6)}+\underbrace{\int_{\R} \rot\vu \cdot (\phi_R^2 \vw) dx}_{(7)}.\notag
\end{align}
Now, we study each term of the right hand side of the expression above.
\begin{itemize}
%%%%%%%%%%%%%%%%%%%%%%%%%%%%%%%%%%%%%%%%%%%%%%%%%%%
\item For the first term above, we write
$$\int_{\R}\Delta(\phi_R^2) (|\vu|^2+|\vw|^2) dx\leq C\int_{\R}\big(|\vn (\phi_R)|^2+\phi_R\Delta(\phi_R)\big) (|\vu|^2+|\vw|^2) dx,$$
and we remark that, by construction, we have 
\begin{equation}\label{Definition_Support_FonctionPhi}
sup(\vn (\phi_R))=sup(\Delta (\phi_R))\subset\{x\in \R: R\leq |x|\leq 2R\}:=\mathcal{C}_R.
\end{equation}
We then obtain, by applying the H\"older inequality with $1=\frac{2}{3}+\frac{1}{3}$:
\begin{eqnarray*}
\int_{\R}\big(|\vn (\phi_R)|^2+\phi_R\Delta(\phi_R)\big) (|\vu|^2+|\vw|^2) dx&\leq &\|\vn (\phi_R)\|^2_{L^3}(\|\vu\|_{L^6(\mathcal{C}_R)}^2+\|\vw\|_{L^6(\mathcal{C}_R)}^2)\\
&&+ \|\phi_R\|_{L^\infty}\|\Delta(\phi_R)\|_{L^{\frac{3}{2}}}(\|\vu\|_{L^6(\mathcal{C}_R)}^2+\|\vw\|_{L^6(\mathcal{C}_R)}^2).
\end{eqnarray*}
Note that, by homogeneity we have $\|\vn (\phi_R)\|_{L^3}=\|\vn \phi\|_{L^3}$ and $\|\Delta(\phi_R)\|_{L^{\frac{3}{2}}}=\|\Delta\phi\|_{L^{\frac{3}{2}}}$ and since $\|\phi_R\|_{L^\infty}\leq C$, we can write
\begin{equation}\label{Estimation_Liouville_1}
\int_{\R}\Delta(\phi_R^2) (|\vu|^2+|\vw|^2) dx\leq C(\|\vu\|_{L^6(\mathcal{C}_R)}^2+\|\vw\|_{L^6(\mathcal{C}_R)}^2).
\end{equation}
%%%%%%%%%%%%%%%%%%%%%%%%%%%%%%%%%%%%%%%%%%%%%%%%%%%
\item For quantity (2) in (\ref{LocalEstimateBothVariables}) we have by the Cauchy-Schwarz inequality and by the H\"older inequality with $\frac{1}{2}=\frac{1}{6}+\frac{1}{3}$:
\begin{eqnarray*}
\int_{\R}\grad (\phi_R^2)\cdot \vw div (\vw) dx&= &2\int_{\R}( \phi_R \grad \phi_R) \cdot \vw div (\vw) dx\leq C\|div(\vw)\phi_R\|_{L^2}\|\vw\cdot \vn \phi_R\|_{L^2}\\
&\leq & C\|\phi_R\|_{L^\infty}\|div(\vw)\|_{L^2}\|\vw\|_{L^6(\mathcal{C}_R)}\|\vn \phi_R\|_{L^3},
\end{eqnarray*}
since $\|\phi_R\|_{L^\infty}\leq C$, $\|\vn \phi_R\|_{L^3}=\|\vn \phi\|_{L^3}\leq C$ by homogeneity and 
$\|div(\vw)\|_{L^2}\leq \|\vw\|_{\dot{H}^1}$, we obtain
\begin{equation}\label{Estimation_Liouville_2}
\int_{\R}\grad (\phi_R^2)\cdot \vw div (\vw) dx\leq C\|\vw\|_{\dot{H}^1}\|\vw\|_{L^6(\mathcal{C}_R)}.
\end{equation}
%%%%%%%%%%%%%%%%%%%%%%%%%%%%%%%%%%%%%%%%%%%%%%%%%%%
\item The term (3) in (\ref{LocalEstimateBothVariables}) is treated as follows:
\begin{eqnarray*} 
\int_{\R} [(\vu \cdot \vec{\nabla})\vu]\cdot (\phi_R^2\vu) dx &=& \sum_{i,j=1}^{3} \int_{B_R} u_j (\partial_{x_j} u_i) ( \phi_R^2 u_i ) dx =  \sum_{i,j=1}^{3} \int_{\R} \phi_R^2 u_j  (\partial_{x_j} u_i) u_i dx\\
&=&  \sum_{i,j=1}^{3} \int_{\R} \phi_R^2 u_j  (\partial_{x_j} \left(\frac{u^{2}_{i}}{2}\right)) dx,
\end{eqnarray*}
and by an integration by parts we obtain
$$\sum_{i,j=1}^{3} \int_{\R} \phi_R^2 u_j  (\partial_{x_j} \left(\frac{u^{2}_{i}}{2}\right)) dx =  - \int_{\R} \vec{\nabla} (\phi_R^2) \cdot \left( \frac{\| \vu\|^2}{2} \vu\right)dx,$$
where we used the fact that $div(\vu)=0$. We thus have
$$\int_{\R} [(\vu \cdot \vec{\nabla})\vu]\cdot (\phi_R^2\vu) dx \leq C\int_{\R} |\vec{\nabla} (\phi_R^2)||\vu|^3dx\leq C\|\phi_R\|_{L^\infty}\int_{\R} |\vec{\nabla} \phi_R||\vu|^3dx.$$
Now, by the H\"older inequalities (with $\frac{1}{\ell}+\frac{3}{q}=1$), by the support properties of $\vn \phi_R$ (see (\ref{Definition_Support_FonctionPhi}) above) and since $\|\phi_R\|_{L^\infty}\leq C$, we have
\begin{equation}\label{eq-aux1} 
\int_{\R} [(\vu \cdot \vec{\nabla})\vu]\cdot (\phi_R^2\vu) dx\leq  \|\vec{\nabla} \phi_R \|_{L^\ell}\| \vu\|_{L^q(\mathcal{C}_R)}^3,
\end{equation}   
moreover, by homogeneity we obtain $\|\vec{\nabla} \phi_R \|_{L^\ell} \leq R^{\frac{3}{\ell}-1}\| \vec{\nabla}\phi \|_{L^\ell}$, and since $\frac{1}{\ell}=1-\frac{3}{q}$, we have $\frac{3}{\ell}-1=2-\frac{9}{q}$ and it follows that $\|\vec{\nabla} \phi_R \|_{L^\ell} \leq R^{2-\frac{9}{q}}\| \vn\phi \|_{L^\ell}$. But, as we have $3\leq q\leq \frac{9}{2}$ by hypotheses, we do have $-1\leq 2-\frac{9}{q}\leq 0$, and since $R>1$ we thus obtain that $R^{2-\frac{9}{q}}\leq 1$ and we can write
 \begin{equation}\label{estim-aux1}
 \|\vec{\nabla} \phi_R \|_{L^\ell}\leq \| \vec{\nabla} \phi \|_{L^\ell}<+\infty. 
 \end{equation}
Thus, coming back to (\ref{eq-aux1}) we have
\begin{equation}\label{Estimation_Liouville_3}
\int_{\R} [(\vu \cdot \vec{\nabla})\vu]\cdot (\phi_R^2\vu) dx \leq  C\| \vu \|^{3}_{L^q(\mathcal{C}R)}.
\end{equation}
%%%%%%%%%%%%%%%%%%%%%%%%%%%%%%%%%%%%%%%%%%%%%%%%%%%
\item The quantity (4) in (\ref{LocalEstimateBothVariables}) can not be treated as previously since the vector field $\vw$ is not divergence free. However, since $div(\vu)=0$ we can still write by an integration by parts:
$$\int_{\R}[(\vu\cdot \grad)\vw] \cdot (\phi_R^2 \vw) dx=-\int_{\R}[(\vu\cdot \grad)(\phi_R^2 \vw)] \cdot  \vw dx,$$
from which we deduce the identities
\begin{eqnarray*}
\int_{\R}[(\vu\cdot \grad)\vw] \cdot (\phi_R^2 \vw) dx&=&-\int_{\R}\big[\big(\vu\phi_R^2)\cdot \grad\big) \vw\big] \cdot  \vw dx-\int_{\R}[\vu\cdot \vn (\phi_R^2)]|\vw|^2dx\\
&=&-\int_{\R}\big[(\vu\cdot \grad) \vw\big] \cdot  (\phi_R^2\vw) dx-\int_{\R}[\vu\cdot \vn (\phi_R^2)]|\vw|^2dx,
\end{eqnarray*}
and we finally obtain
$$\int_{\R}[(\vu\cdot \grad)\vw] \cdot (\phi_R^2 \vw) dx=-\frac{1}{2}\int_{\R}[\vu\cdot \vn (\phi_R^2)]|\vw|^2dx.$$
We can now write, by the Young inequalities for the sum:
\begin{eqnarray*}
\int_{\R}[(\vu\cdot \grad)\vw] \cdot (\phi_R^2 \vw) dx&\leq &\sum_{j=1}^3\int_{\R}|u_j \phi_R (\partial_{x_j}\phi_R)||\vw|^2dx\leq \sum_{j=1}^3\int_{\R}\big(4|u_j\partial_{x_j}\phi_R|^2+\frac{|\phi_R|^2}{4} \big)|\vw|^2dx\\
&\leq &C\int_{\R}|u|^2|\vn\phi_R|^2|\vw|^2 dx+\frac{3}{4}\int_{\R}|\phi_R\vw|^2dx.
\end{eqnarray*}
By the H\"older inequalities with $1=\frac{1}{3}+\frac{1}{3}+\frac{1}{3}$ we have
$$\int_{\R}[(\vu\cdot \grad)\vw] \cdot (\phi_R^2 \vw) dx\leq C\|\vu\|_{L^6}^2\|\vn\phi_R\|_{L^6}^2\|\vw\|_{L^6}^2+\frac{3}{4}\|\phi_R\vw\|_{L^2}^2,$$
and since we have, by homogeneity, that $\|\vn\phi_R\|_{L^6}^2=R^{-1}\|\vn\phi\|_{L^6}^2\leq CR^{-1}$, we write
\begin{equation}\label{Estimation_Liouville_4}
\int_{\R}[(\vu\cdot \grad)\vw] \cdot (\phi_R^2 \vw) dx\leq CR^{-1}\|\vu\|_{L^6}^2\|\vw\|_{L^6}^2+\frac{3}{4}\|\phi_R\vw\|_{L^2}^2.
\end{equation}
%%%%%%%%%%%%%%%%%%%%%%%%%%%%%%%%%%%%%%%%%%%%%%%%%%%
\item For the term (5) in (\ref{LocalEstimateBothVariables}), using the divergence free property for $\vu$ we have
\begin{eqnarray} \nonumber
\int_{\R} \vec{\nabla}p\cdot ( \phi_R^2 \vu )dx&=& \sum_{i=1}^{3}\int_{\R}(\partial_{x_i} p) \phi_R^2 u_i dx=- \sum_{i=1}^{3}\int_{\R} p \partial_{x_i} (\phi_R^2 u_i)dx \\
&=& - \sum_{i=1}^{3}\int_{\R} p \partial_{x_i} (\phi_R^2)  ( u_i) dx = - \int_{\R}\vec{\nabla}(\phi_R^2)\cdot(p \vu) dx. \label{ipp4}
\end{eqnarray}
Using again the H\"older inequalities (with $\frac{1}{\ell}+\frac{3}{q}=1$), by the estimate (\ref{estim-aux1})  -as we are working with the same set of exponents- and taking into account the support properties of $\vn\phi_R$, we obtain
\begin{eqnarray}
&&\int_{\R} \vec{\nabla}p\cdot ( \phi_R^2 \vu )dx\leq  C\|\phi_R\|_{L^\infty}\int_{\R}| \vec{\nabla} \phi_R| |p| |\vu| dx \notag\\ 
&&\leq  C\|\phi_R\|_{L^\infty}\|\vec{\nabla} \phi_R \|_{L^\ell} \left( \int_{\mathcal{C}_R} (|p| |\vu|)^{\frac{q}{3}} dx \right)^{\frac{3}{q}}\leq  C \left( \int_{\mathcal{C}_R} (|p| |\vu|)^{\frac{q}{3}} dx \right)^{\frac{3}{q}}.\label{estim-aux2}
\end{eqnarray}
Now, since by hypothesis we have $\vu\in L^q(\R)$ and since by the estimate (\ref{Estimation_Pression_Lq2}) we have $p\in L^{\frac{q}{2}}(\R)$, thus by the H\"older inequalities (with $\frac{2}{q}+\frac{1}{q}=\frac{3}{q}$) we can write
\begin{equation}\label{Estimation_Liouville_5}
\int_{\R} \vec{\nabla}p\cdot ( \phi_R^2 \vu )dx\leq C \|p\|_{L^\frac{q}{2}(\mathcal{C}_R)} \|\vu\|_{L^q(\mathcal{C}_R)}.
\end{equation}
%%%%%%%%%%%%%%%%%%%%%%%%%%%%%%%%%%%%%%%%%%%%%%%%%%%
\item For the  term $(6)$ in (\ref{LocalEstimateBothVariables}), by the Cauchy-Schwarz inequality and the Young inequalities for the sum, we obtain
\begin{eqnarray*}
\frac{1}{2}\int_{\R} \vw\cdot[\vn(\phi_R^2)\wedge \vu] dx&\le &  \int_{\R}|\phi_R| |\grad(\phi_R)||\vw||\vu| dx  \le \|\phi_R\vw\|_{L^2}\|\grad(\phi_R)\vu\|_{L^2}\\
&\le &\frac{\|\phi_R\vw\|_{L^2}^2}{4}+4\|\grad(\phi_R)\vu \|_{L^2}^2,
\end{eqnarray*}
and applying the Hölder inequality with $\frac{1}{2}=\frac{1}{3}+\frac{1}{6}$ in the last term above, we have
$$\leq \frac{\|\phi_R\vw\|_{L^2}^2}{4}+ 4\|\grad(\phi_R) \|_{L^3}^2\|\vu\|_{L^6(\mathcal{C}_R)}^2.$$
Recalling that, by homogeneity we have $\|\vn (\phi_R)\|_{L^3}=\|\vn \phi\|_{L^3}<+\infty,$ we finally obtain
\begin{equation}\label{Estimation_Liouville_6}
\frac{1}{2}\int_{\R}\rot \vw\cdot( \phi_R^2 \vu) dx\le 
\frac{\|\phi_R\vw\|_{L^2}^2}{4}+ C \|\vu\|_{L^6(\mathcal{C}_R)}^2.\\[5mm]
\end{equation}
\item For the last term of  (\ref{LocalEstimateBothVariables}) by  the Cauchy-Schwartz inequality we have
\begin{equation*}
\int _{\R} \rot\vu \cdot (\phi_R^2 \vw) dx\le \|\phi_R\rot\vu\|_{L^2}\|\phi_R\vw\|_{L^2} 
\le 2\|\phi_R\grad \otimes \vu\|_{L^2}\|\phi_R\vw\|_{L^2}. 
\end{equation*}
Hence, by the Young inequality for the sum (introducing a small technical parameter $0<\mathfrak{e}<1$) we observe that
\begin{equation}\label{Estimation_Liouville_7}
\int _{\R} \rot\vu \cdot (\phi_R^2 \vw) dx\le \mathfrak{e}\|\phi_R\grad\otimes\vu\|_{L^2}^2+C(\mathfrak{e})\|\phi_R\vw\|_{L^2}^2.
\end{equation}

\end{itemize}{

%%%%%%%%%%%%%%%%%%%%%%%%%%%%%%%%%%%%%%%%%%%%%%%%%%%
With all these estimates (\ref{Estimation_Liouville_1})-(\ref{Estimation_Liouville_7}) at hand, we come back to (\ref{LocalEstimateBothVariables}) and we have
\begin{align*}
\int_{\R} \phi_R^2| \grad\otimes \vu|^2dx&+\int_{\R} \phi_R^2| \grad\otimes \vw|^2dx+\int_{\R} \phi_R^2| div( \vw)|^2 dx+\kappa \int_{\R} \phi_R^2|\vw|^2dx\qquad \notag\\
&\leq C(\|\vu\|_{L^6(\mathcal{C}_R)}^2+\|\vw\|_{L^6(\mathcal{C}_R)}^2)+ C\|\vw\|_{\dot{H}^1}\|\vw\|_{L^6(\mathcal{C}_R)}+C\| \vu \|^{3}_{L^q(\mathcal{C}_R)}\notag\\
&+ CR^{-1}\|\vu\|_{L^6}^2\|\vw\|_{L^6}^2+\frac{3}{4}\|\phi_R\vw\|_{L^2}^2+ C \|\vu\|_{L^\frac{q}{2}(\mathcal{C}_R)} \|p\|_{L^q(\mathcal{C}_R)}+\frac{1}{4}\|\phi_R\vw\|_{L^2}^2+ C \|\vu\|_{L^6(\mathcal{C}_R)}^2\\
&+\mathfrak{e}\|\phi_R \grad\otimes\vu\|_{L^2}^2+  C(\mathfrak{e})\|\phi_R \vw\|_{L^2}^2.\notag
\end{align*}

We now remark that the quantity $\displaystyle{\int_{\R} \phi_R^2| div( \vw)|^2 dx}$ in the first line above is positive and that the quantities $\|\phi_R \grad\otimes\vu\|_{L^2}^2$ and $\|\phi_R\vw\|_{L^2}$ are bounded (recall moreover that $\kappa\gg 1$ and $0<\mathfrak{e}<1$), thus we can write 
\begin{align}
(1-\mathfrak{e})\int_{\R} \phi_R^2| \grad\otimes \vu|^2dx&+\int_{\R} \phi_R^2| \grad\otimes \vw|^2dx+(\kappa-C(\mathfrak{e})-1) \|\phi_R\vw\|_{L^2}^2\qquad \notag\\
&\leq C(\|\vu\|_{L^6(\mathcal{C}_R)}^2+\|\vw\|_{L^6(\mathcal{C}_R)}^2)+ C\|\vw\|_{\dot{H}^1}\|\vw\|_{L^6(\mathcal{C}_R)}+C\| \vu \|^{3}_{L^q(\mathcal{C}_R)}\label{Derniere_Estimation}\\
&+CR^{-1}\|\vu\|_{L^6}^2\|\vw\|_{L^6}^2+C \|p\|_{L^\frac{q}{2}(\mathcal{C}_R)} \|\vu\|_{L^q(\mathcal{C}_R)}+ C \|\vu\|_{L^6(\mathcal{C}_R)}^2.\notag
\end{align}
At this point we note that, since $\vu, \vw\in L^6(\R)$, $\vw\in L^2(\R)$ and $\vu\in L^q(\R)$ and $p\in  L^{\frac{q}{2}}(\R)$ with $3\leq q\leq \frac{9}{2}$, then the quantities
$$\|\vu\|_{L^6(\mathcal{C}_R)}, \quad\|\vw\|_{L^6(\mathcal{C}_R)},  \quad\|\vu\|_{L^q(\mathcal{C}_R)}, \quad R^{-1}\|\vu\|_{L^6}^2\|\vw\|_{L^6}^2, \quad \|p\|_{L^\frac{q}{2}(\mathcal{C}_R)},$$
 which are present in each term in the right-hand side of the previous estimate, will tend to $0$ if $R\to +\infty$ and then all the right-hand side of (\ref{Derniere_Estimation}) above will tend to $0$. For the left-hand side of (\ref{Derniere_Estimation}), if we let $R\to +\infty$ we obtain the quantities $\|\vu\|_{\dot{H}^1}$, $\|\vw\|_{\dot{H}^1}$ and $\|\vw\|_{L^2}$ which will thus tend to $0$. Since we have the Sobolev embedding $\dot{H}^{1}(\R)\subset L^6(\R)$, we easily deduce that $\vu=0$ and $\vw=0$. 
The proof of Theorem \ref{Theorem_Liouville_type} is finished. \hfill $\blacksquare$
%%%%%%%%%%%%%%%%%%%%%%%%%%%%%%%%%%%%%%%%%%%%%%%%%%%

%%%%%%%%%%%%%%%%%%%%%%%%%%%%%%%%%%%%%%%%%%%%%%%%%%%
\end{document}